\documentclass[12pt]{amsart}

\usepackage{amsfonts}
\usepackage{amssymb}
\usepackage{amsmath}
\usepackage{amsthm}
\usepackage{latexsym}
\usepackage{amsthm}
\usepackage{graphicx}
\usepackage{amsmath}
\usepackage{amsfonts}
\usepackage{tikz}
\usetikzlibrary{matrix,arrows}
\usepackage{amscd}
\usepackage{amssymb}
\usepackage{xy}
\usepackage{units}
\usepackage{fancyhdr}

\def\aint{\mathop{-\kern-10.0pt\int}}

\newtheorem{theorem}{Theorem}[subsection]

\newtheorem{definition}[theorem]{Definition}

\newtheorem{lemma}[theorem]{Lemma}

\newtheorem{proposition}[theorem]{Proposition}
\newtheorem{remark}[theorem]{Remark}

\numberwithin{equation}{subsection}

\begin{document}

\title{Scalar invariants of surfaces in conformal 3-sphere via Minkowski spacetime}
\author{Jie Qing, Changping Wang, and Jingyang Zhong}

\begin{abstract} For a surface in 3-sphere, by identifying the conformal round 3-sphere as the projectivized positive light cone in Minkowski 5-spacetime, 
we use the conformal Gauss map and the conformal transform to construct the associate homogeneous 4-surface in Minkowski 5-spacetime. We then 
derive the local fundamental theorem for a surface in conformal round 3-sphere from that of the associate 4-surface in Minkowski 5-spacetime. 
More importantly, following the idea of Fefferman and Graham \cite{FG-1, FG-2}, we construct local scalar invariants for a surface in conformal round 3-sphere. 
One distinct feature of our construction is to link the classic work of Blaschke \cite{blaschke}, Bryan \cite{bryant} and Fefferman-Graham \cite{FG-1, FG-2}.

\end{abstract}
\maketitle
%\tableofcontents

\section{Introduction}

It is well-known that all local scalar invariants of a (pseudo-)Riemannian metric are Weyl invariants, based on Weyl's classical invariant theory for the orthogonal groups. 
A conformal structure on a manifold is described by an equivalent class of conformal Riemannian metrics. Two metrics $g_1$ and $g_2$ on a manifold $\textup{M}$ 
are conformal to each other if  $g_1 = \lambda^2 g_2$ for some positive smooth function $\lambda$ on $\textup{M}$. There are several ways to set the theory of local conformal invariants, but it is no longer straightforward to account for local scalar conformal invariants because of the lack of Weyl Theorem for the group of conformal transformations. To tackle such problem, in the seminal paper \cite{FG-1} in 1980's,  Fefferman and Graham described the ingenious construction of a Ricci-flat homogeneous Lorentzian ambient spacetime for a given conformal manifold, where the conformal manifold is represented by the homogeneous null hypersurface in the ambient spacetime. Their construction was motivated by the model case in which  the conformal round sphere $\mathbb{S}^n$ is the projectivized positive light cone $\mathbb{N}^{n+1}_+$ in Minkowski spacetime $\mathbb{R}^{1,n+1}$. In \cite{FG-1}, Fefferman and Graham initiated the program to use local scalar (pseudo-)Riemannian invariants of the ambient metrics at the homogeneous null hypersurface to fully account for local scalar conformal invariants. Readers are referred to their recent expository paper \cite{FG-2} to learn all the developments of this program  (cf. also, \cite{BEG, gover}). This program also has lead to many significant advances in the global theory of  conformal geometry, particularly via conformally invariant PDEs. \\

In this paper we want to build the model case to the study of local scalar invariants of submanifolds in a conformal manifold 
in the way that follows the approach in \cite{FG-1}. The model case for us is
to study 2-surfaces $\hat x$ in the conformal round 3-sphere $(\mathbb{S}^3, [g_0])$. As in \cite{FG-1}, the conformal round 3-sphere is represented by the positive 
light cone $\mathbb{N}^4_+$ in Minkowski 5-spacetime $\mathbb{R}^{1,4}$.  Given an immersed  surface
$$
\hat x: \textup{M}^2\to\mathbb{S}^3
$$
or equivalently
$$
y = (1, \hat x): \textup{M}^2\to \mathbb{N}^4_+,
$$
to incorporate all metrics in $[g_0]$ on 3-sphere we consider the homogeneous extension
$$
x^{\mathbb{N}} = \alpha (1, \hat x): \mathbb{R}^+\times\textup{M}^2\to \mathbb{N}^4_+\subset\mathbb{R}^{1,4}.
$$
Then we will use the conformal Gauss map $\xi$ of $\hat x$ to choose a canonical null vector $y^*$ at each given point $y \in x^{\mathbb{N}}\subset \mathbb{N}^4_+$
to extend $x^{\mathbb{N}}$ further into a homogeneous timelike 4-surface
$$
\tilde x=\alpha y + \alpha\rho y^*: \mathbb{R}^+\times\mathbb{R}^+\times\textup{M}^2\to\mathbb{R}^{1,4}.
$$
We will also consider the associate ruled 3-surface 
$$
x^+ = \frac 1{\sqrt 2} (e^ty + e^{-t}y^*): \mathbb{R}\times\textup{M}^2\to \mathbb{H}^4\subset\mathbb{R}^{1,4}
$$
where $\mathbb{H}^4$ is the hyperboloid in Minkowski 5-spacetime. The main idea, inspired by the work \cite{FG-1, FG-2}, is to use the geometry of the associate
4-surface $\tilde x$ in Minkowski spacetime $\mathbb{R}^{1,4}$ (the associate ruled 3-surface $x^+$ in the hyperboloid $\mathbb{H}^4$ and the spacelike surface as the
image of the conformal Gauss map $\xi$ in the de Sitter spacetime $S^{3,1}$ in Minkowski spacetime $\mathbb{R}^{1,4}$) to study 
the geometry of the surface $\hat x$ in the conformal round 3-sphere $\mathbb{S}^3$.\\

Our approach facilitates proofs of the local fundamental theorems (cf. Theorem \ref{local-fundamental} and \cite{cpw-1, cpw-2}) 
and produces local scalar invariants of surfaces in the conformal round 3-sphere. 
The second is more interesting and helpful to find appropriate PDE problems to study the surfaces.
The study of Willmore surfaces indeed exemplifies well that how important and central those problems are in the theory of surfaces in general \cite{blaschke, bryant, LY, MN}.  \\

We should remark that the key to our construction of associate surfaces is the conformal Gauss map $\xi$ to a given surface $\hat x$ in the conformal round 3-sphere. 
The conformal Gauss maps have been introduced in several contexts (cf. \cite{blaschke, bryant, rigoli}). We are searching for a definition that fits into the context of ambient spaces
of Fefferman and Graham (cf. Lemma \ref{xi-1} and Lemma \ref{xi-2}). It is fascinating to see how Blaschke \cite{blaschke} introduced the conformal Gauss map as the 
map representing the family of mean curvature 2-spheres of the surface $\hat x$ and the conformal transform $\hat x^*$ (cf. Definition \ref{def-conf-trans}) as the other envelope surface
of the conformal Gauss map. One technical assumption for the null vector $y^*$ to be well defined at each point $y\in x^{\mathbb{N}}$ is to require that the conformal Gauss map of the 
surface $\hat x$ induces a spacelike surface in the de Sitter spacetime $\mathbb{S}^{1,3}$, which is equivalent to that the surface $\hat x$ is free of 
umbilical point in the conformal 3-sphere $\mathbb{S}^3$.\\

It is nice to know that in our construction the associate 4-surface $\tilde x$ in Minkowski spacetime $\mathbb{R}^{1,4}$ is a minimal 4-surface (of vanishing mean curvature) 
if and only if the 2-surface $\hat x$ is a Willmore surface with no umbilical point in $\mathbb{S}^3$ (cf. Theorem \ref{Th:tilde-H}). The same statement also holds for the associate 
ruled 3-surface $x^+$ in the hyperboloid $\mathbb{H}^4$ (cf. Theorem \ref{Th:+-H}) as well as the conformal Gauss map surface $\xi$ in de Sitter sapcetime $\mathbb{S}^{1,3}$ 
(cf. Theorem \ref{conf-gauss-surface}).\\

Upon realizing that a different representative $\lambda^2 g_0$ in the conformal class $[g_0]$ on $\mathbb{S}^3$ is equivalent to a different
parametrization for the associate surface
\begin{equation}\label{general-param}
\tilde x = \alpha y_\lambda + \alpha\rho y^*_\lambda:\mathbb{R}^+\times\mathbb{R}^+\times\textup{M}^2\to\mathbb{R}^{1,4}, 
\end{equation}
where $y_\lambda = \hat \lambda (1, \hat x)$ and $\hat \lambda = \lambda\circ\hat x$ for a conformal factor $\lambda$, the real issue is how we use the geometry of the surface 
$\hat x$ in the 3-sphere $(\mathbb{S}^3, \lambda^2g_0)$ to calculate the geometry of the associate surface $\tilde x$. The solution is to use the following 3-sphere 
$\mathbb{S}^3_\lambda$ in the positive light cone $\mathbb{N}^4_+$:
\begin{equation}
\lambda (1, x): \mathbb{S}^3\to \mathbb{N}^4_+
\end{equation}
as the realization of $(\mathbb{S}^3, \lambda^2g_0)$.
For the convenience of readers we present the calculations of the geometry of $\mathbb{S}^3_\lambda$ as a spacelike 3-surface in Minkowski spacetime  
in the Appendix \ref{gauss}. But it starts with the following observation.

\begin{lemma} Suppose that $\hat x: \textup{M}^2\to \mathbb{S}^3$ is an immersed surface and $\lambda^2g_0$ is a conformal metric in the round conformal class $[g_0]$
on $\mathbb{S}^3$. Then 
\begin{equation}
\xi = H_\lambda y_\lambda + \overset{\rightarrow}{\bf n}_\lambda,
\end{equation}
where $H_\lambda$ is the mean curvature of $\hat x$ in $(\mathbb{S}^3, \lambda^2g_0)$ 
and $\overset{\rightarrow}{\bf n}_\lambda$ is the unit normal to $y_\lambda$ in $\mathbb{S}^3_\lambda\subset\mathbb{N}^4_+$.
\end{lemma}

Using the calculations in Appendix \ref{gauss}, we are able to show in the proof of Theorem \ref{main-inv} that the data $\{m, \omega^\lambda, \Omega_\lambda, 
\Omega^*_\lambda\}$ that determine the first and second fundamental forms of the associate surface $\tilde x$ in Minkowski spacetime $\mathbb{R}^{1, 4}$ 
can all be expressed in terms of covariant derivatives of
the curvature of the surface $\hat x$ in $(\mathbb{S}^3, \lambda^2g_0)$ and the covariant derivatives of curvature of $(\mathbb{S}^3, \lambda^2g_0)$ (including 0th order).  
In the exact same spirit as in Fefferman and Graham \cite{FG-1, FG-2},  our construction of associate surfaces $\tilde x$ provides a way to capture local scalar conformal invariants of 
a surface $\hat x$ . Namely, one can obtain local scalar conformal invariants of the surface $\hat x$ in the conformal round 3-sphere by computing the 
local scalar (pseudo-)Riemannian invariants of the associate surface $\tilde x$ at the homogeneous surface $x^{\mathbb{N}}$ in the light cone in Minkowski 5-spacetime. 
The first non-trivial one is 
\begin{equation}\label{laplacian-tilde}
\tilde\Delta\tilde H|_{\rho=0}  = 2\alpha^{-3} (\Delta_\lambda H_\lambda  + |\overset{\circ}{II}_\lambda|^2H_\lambda + (\overset{\circ}{II}_\lambda)^{ij}(R^\lambda)_{i3j3} 
-  (R^\lambda)_{3i,}^{\quad i} )
\end{equation} 
in a general parametrization \eqref{general-param},  where $(R^\lambda)_{i3j3}$ and $(R^\lambda)_{3i}$ are the Riemann curvature and Ricci curvature of the metric $\lambda^2g_0$ 
on $\mathbb{S}^3$. Due to the homogeneity of $\tilde x$ we automatically have
\begin{equation}\label{intr-willmore}
{\mathcal H}_\lambda = \Delta_\lambda H_\lambda  + |\overset{\circ}{II}_\lambda|^2H_\lambda +
(\overset{\circ}{II}_\lambda)^{ij}(R^\lambda)_{i3j3}  -  (R^\lambda)_{3i,}^{\ \ i}= \hat\lambda^{-3} ( \Delta H + |\overset{\circ}{II}|^2H)
\end{equation}
which is the curvature that vanishes if and only if the surface $\hat x$ is Willmore. Notice that extra curvature terms do not show up when we work with either the round metric $g_0$ 
or the Euclidean metric. Similar formulas have appeared in the literature \cite{HL, GR, gover-s}. \\

We also calculate in Section \ref{calculations} some other conformal scalar invariants of higher orders:
\begin{equation}\label{intro-norm-co-der}
\aligned 
|\nabla\tilde h|^2|_{\rho =0} & =\alpha^{-4}( |\nabla\Omega_\lambda|^2 + 8|dH_\lambda|^2 + 2 Ric^\lambda(\overset{\rightarrow}{\bf n}_\lambda, \nabla H_\lambda) + 3 H^2_\lambda|\Omega_\lambda|^2 \\ & \quad + 3 K^T_\lambda |\Omega_\lambda|^2+ 6 \Omega_\lambda\cdot \text{Hess} (H_\lambda)) \endaligned
\end{equation}
(cf. \eqref{norm-co-der}, where $K^T_\lambda$ is the sectional curvature of $(\mathbb{S}^3, \lambda^2 g_0)$ at the tangent plane to the surface $\hat x$,  and 
\begin{equation}\label{intro-double-laplace}
\aligned
\tilde\Delta & \tilde\Delta\tilde H|_{\rho=0} = 8\alpha^{-5} (\Delta_\lambda{\mathcal H}_\lambda +9 |\omega^\lambda|^2{\mathcal H}_\lambda - 3 \text{Div}({\omega^\lambda})
{\mathcal H}_\lambda \\ & \quad - 6\omega^\lambda(\nabla{\mathcal H}_\lambda) 
- 6 {\mathcal H}_\lambda| \overset{\circ}{II}_\lambda|^{-2}\overset{\circ}{II}_\lambda\cdot\Omega^*_\lambda),\endaligned
\end{equation}
where $\omega^\lambda = <dy_\lambda, y^*_\lambda>$ and $\Omega^*_\lambda = - <dy^*_\lambda, d\xi>$ are parts of the data that determine the geometry of the associate surface $\tilde x$
and are given in \eqref{omega} and \eqref{Omega-star-ij} as invariants of the surface $\hat x$ in $(\mathbb{S}^3, \lambda^2g_0)$.
\\

To end the introduction we remark that, for the sake of the production of local scalar invariants, the assumption of having no umbilical point in our construction 
is not an issue.

\section{The associate surfaces in $\mathbb{R}^{1,4}$}\label{Sect:2}

In this section we introduce the associate surfaces in Minkowski space $\mathbb{R}^{1,4}$ for a given surface
$\hat x: \textup{M}^2\to \mathbb{S}^3$.  
We then show that such associate surface is canonical in doing conformal geometry for the surface $\hat x$. The construction relies on the conformal Gauss map 
and the conformal transform of $\hat x$. It is also very interesting to see how Blaschke and Bryant came to the conformal Gauss map and the conformal transform in very different 
perspectives \cite{blaschke, bryant}.

\subsection{Surfaces in 3-sphere}

Suppose that 
$$
\hat x: \textup{M}^2\to \mathbb{S}^3\subset \mathbb{R}^4
$$
is an immersed surface with isothermal coordinate $(u^1, u^2)$. Let 
$$
{\bf n}: \textup{M}^2\to \mathbb{R}^4
$$
be the unit normal vector at each point on the surface. Then we obtain the first fundamental form 
\begin{equation}\label{x-hat-I}
I = <d\hat x, d\hat x> = E|du|^2
\end{equation}
and the second fundamental form 
\begin{equation}\label{x-hat-II}
II = - <d\hat x, d{\bf n}>  = e(du^1)^2 + 2f du^1du^2 + g(du^2)^2.
\end{equation}
Hence the mean curvature of the surface in 3-sphere is 
\begin{equation}\label{x-hat-H}
H = \frac 1{2E}(e + g)
\end{equation}
and the Gaussian curvature of the surface is  
\begin{equation}\label{x-hat-guass}
K = \frac {eg - f^2}{E^2} + 1.
\end{equation}
Notice that 
\begin{equation}\label{x-hat-gauss-map}
\left\{\aligned {\bf n}_{u^1} & = - \frac eE \hat x_{u^1} - \frac fE \hat x_{u^2}\\
                        {\bf n}_{u^2} & = - \frac fE \hat x_{u^1} - \frac gE \hat x_{u^2}.\endaligned
\right.
\end{equation}            
\\
 
If one takes another conformal metric $\lambda^2 g_0$ on the 3-sphere $\mathbb{S}^3$, where $\lambda$ is a positive function on $\mathbb{S}^3$, then
the first fundamental form for the surface $\hat x$ is
\begin{equation}\label{lambda-1}
I_\lambda = \hat\lambda^2 I, 
\end{equation}
where $\hat\lambda = \lambda\circ\hat x$ and the second fundamental form is
\begin{equation}\label{lambda-2}
II_\lambda = \hat\lambda II - \lambda_{\bf n} I,
\end{equation}
where $\lambda_{\bf n} = {\bf n}(\lambda)$. Hence 
\begin{equation}\label{after-lambda}
H_\lambda = \hat\lambda^{-1} (H - \frac {\lambda_{\bf n}}{\hat\lambda}) \text{ and } \overset{\circ} II_\lambda = \hat\lambda \overset{\circ}{II},
\end{equation}
where $\overset{\circ}{II}$ is the traceless part of the second fundamental form $II$. Here we see the easy scalar conformal invariant
$|\overset{\circ} II|^{2}$, which can be considered to be the counter part of the square of the length of Weyl curvature on a conformal manifold.

\subsection{Minkowski 5-spacetime}

Let $\mathbb{R}^{1,4}$ be the Minkowski 5-spacetime, where we use the notation
$$
\mathbb{R}^{1,4} = \{(t,  x): t\in\mathbb{R} \text{ and } x\in \mathbb{R}^4\}
$$
with the Lorentz inner product  
$$
< (t, x), (s, y)> = -st + x\cdot  y.
$$
Recall the positive light cone is given by
$$
\mathbb{N}^4_+ = \{(t,  x)\in \mathbb{R}^{1,4}: -t^2 + |x|^2 = 0 \text{ and } t > 0\};
$$
the hyperboloid is given as
$$
\mathbb{H}^4 = \{(t, x)\in \mathbb{R}^{1,4}: -t^2 +|x|^2 = -1 \text{ and } t>0\};
$$
and the de Sitter 4-spacetime is given as
$$
\mathbb{S}^{1,3} = \{(t, x)\in \mathbb{R}^{1,4}: -t^2 + |x|^2 = 1\}.
$$
\\

Given a surface $\hat x: \textup{M}^2\to\mathbb{S}^3\subset\mathbb{R}^4$, we may consider the 2-surface 
$$
y = (1, \hat x): \textup{M}^2\to \mathbb{N}^4_+\subset \mathbb{R}^{1,4}
$$
and the homogeneous extension
$$
x^{\mathbb{N}} = \alpha y: \mathbb{R}^+\times\textup{M}^2\to \mathbb{N}^4_+\subset \mathbb{R}^{1,4}
$$
for $\alpha\in\mathbb{R}^+$.  There does not seem to be a way of doing ``geometry" of the homogeneous 3-surface $x^{\mathbb{N}}$ in the positive light cone
$\mathbb{N}^4_+$. \\

To motivate our choice of the associate surface in $\mathbb{R}^{1,4}$ of $\hat x$ we first introduce the so-called homogeneous
coordinate for $\mathbb{R}^{1,4}$
used in the ambient space construction of Fefferman and Graham \cite{FG-1, FG-2}, that is, 
\begin{equation}\label{f-g-coor}
(t, x) = x^0 (1, \hat x) + x^0x^\infty\frac 12(1, -\hat x)
\end{equation}
where
$$
\left\{\aligned x^0 & = \frac 1{2} (r +t)\\
x^0x^\infty & =  (-r+t)\endaligned\right.
$$
and $r = |x|$ and $x = r\hat x$. In this coordinate the Minkowski metric is
$$
\tilde {\mathcal G}_0 = - 2x^\infty (dx^0)^2 - 2x^0dx^0dx^\infty + (x^0)^2(1 - \frac {x^\infty}2)^2 g_0(\hat x).
$$
Hence, given a surface $\hat x: \textup{M}^2\to\mathbb{S}^3$, we are looking to construct an associate homogeneous timelike
4-surface
\begin{equation}\label{assoc-surface}
\tilde x = \alpha  y + \alpha\rho  y^*: \mathbb{R}^+\times\mathbb{R}^+\times\textup{M}^2\to \mathbb{R}^{1, 4} 
\end{equation}
if we can have canonically the null vector  $y^*$ at a given null  position $y$ on $x^{\mathbb{N}}$. It is clear that the associate surface $\tilde x$ is ruled by the positive quadrants of timelike 2-planes in Minkowski spacetime.
One may consider the intersection of $\tilde x$ with the hyperboloid $\mathbb{H}^4$:
\begin{equation}\label{assoc-in-hyper}
x^+ = \frac 1{\sqrt{2}} (e^t y + e^{-t} y^*): \mathbb{R}\times\textup{M}^2\to \mathbb{H}^4,
\end{equation}
which is called the associate ruled 3-surface since it is a 3-surface in hyperbolic 4-space ruled by geodesics lines. 
Recall that a geodesic line in the hyperboloid $\mathbb{H}^4$ is the intersection of 
the hyperboloid with a timelike 2-subspaces in Minkowski spacetime. 
In the following we will introduce the canonical choice of such $y^*$.\\

\subsection{Conformal Gauss maps}

Let us consider any unit spacelike normal vector to the homogeneous null 3-surface $x^{\mathbb{N}} = \alpha y$ 
in $\mathbb{N}^4_+\subset\mathbb{R}^{1,4}$. That is to ask a unit spacelike 5-vector
$\xi$ to satisfy
\begin{equation}\label{enveloping}
<\xi, x^{\mathbb{N}}>  = 0, \quad <\xi, x^{\mathbb{N}}_{u^1}>  = 0, \quad <\xi,  x^{\mathbb{N}}_{u^2}> = 0,
\end{equation}
which implies that
$$
\xi = a y + \overset{\rightarrow}{\bf n},
$$
where $\overset{\rightarrow}{\bf n} = (0, {\bf n})$ is the unit normal to the surface $\hat x$ in the standard unit round 3-sphere in $\{1\}\times\mathbb{R}^4 \subset\mathbb{R}^{1,4}$. 
It turns out that there is a unique choice if we insist that the map
$$
\xi: \textup{M}^2\to \mathbb{S}^{1,3}\subset \mathbb{R}^{1,4}
$$
is (weakly) conformal. Namely we have

\begin{lemma} \label{xi-1} Suppose that $\hat x:\textup{M}^2\to\mathbb{S}^3$ is an immersed surface. Then, 
for a unit normal vector $\xi$ to the homogeneous null 3-surface $x^{\mathbb{N}} = \alpha y: \mathbb{R}^+\times\textup{M}^2\to \mathbb{N}^4_+\subset \mathbb{R}^{1,4}$, 
$$
<\xi_{u^1}, \xi_{u^2}> = 0
$$
if and only if 
$$
\xi = H y + \overset{\rightarrow}{\bf n}
$$
and
\begin{equation}\label{mobius-metric} 
<d\xi, d\xi> = \frac 12 E|\overset{\circ}{II}|^2 |du|^2.
\end{equation}
\end{lemma}

\proof It is simply a straightforward calculation. We know
$$
\xi_{u^i} = a_{u^i}(1, \hat x) + a(0, \hat x_{u^i}) + (0, {\bf n}_{u^i}).
$$
Hence we have
$$
<\xi_{u^1}, \xi_{u^2}> = -2a f + \frac 1E (fe+fg) = 0,
$$
which is equivalent to $a = H$. For the rest we calculate
\begin{equation}\label{mobius-metric-explicit}
<\xi_{u^1}, \xi_{u^1}>  = <\xi_{u^2}, \xi_{u^2}> = \frac 1{E^2}(f^2 +(\frac {e-g}2)^2)E.
\end{equation}
\endproof

Another way to identify a unique unit spacelike normal vector to the homogeneous null 3-surface $x^{\mathbb{N}} = \alpha y: \mathbb{R}^+\times\textup{M}^2\to\mathbb{N}^4_+$
is the following:

\begin{lemma}\label{xi-2} Suppose that $\hat x: \textup{M}^2\to\mathbb{S}^3$ is an immersed surface. Then, for a unit spacelike normal vector $\xi$ to 
$x^{\mathbb{N}} = \alpha y:\mathbb{R}^+\times\textup{M}^2\to \mathbb{N}^4_+\subset \mathbb{R}^{1,4}$,
$$
\xi = H y+ \overset{\rightarrow}{\bf n}
$$
if and only if
\begin{equation}\label{integrability}
< \Delta \xi, y> = 0.
\end{equation}
\end{lemma}

\proof We simply calculate, for $\xi = a(1, \hat x) + (0, {\bf n})$,
$$
\Delta_0\xi = \xi_{u^1 u^1} + \xi_{u^2 u^2} = (\Delta_0 a)(1, \hat x) + 2\nabla a(0, \nabla \hat x) + a(0, \Delta_0\hat x) + (0, \Delta_0 {\bf n})
$$
and
$$
<\Delta_0\xi, (1, \hat x)> = -2aE + 2HE.
$$
Notice that $\Delta = E^{-1}\Delta_0$.
\endproof

Before we give a formal definition of the conformal Gauss map we want to make a remark that \eqref{integrability} is the integrability condition for the unit vector field $\xi$ to be 
the conformal Gauss map (up to a sign) for the surface $\hat x$. This turns out to be the easiest way to see that $\hat x$ is Willmore if and only if 
the conformal Gauss map $\xi$ of $\hat x$ is also the conformal Gauss map (up to a sign) of the conformal transform $\hat x^*$ (cf. Definition \ref{def-conf-trans}). 

\begin{definition} Suppose that $\hat x: \textup{M}^2\to \mathbb{S}^3$ is a surface. Then 
we will call 
\begin{equation}\label{conformal-gauss-map}
\xi =  H y+ \overset{\rightarrow}{\bf n}: \textup{M}^2\to\mathbb{S}^{1,3}\subset\mathbb{R}^{1,4}
\end{equation}
the conformal Gauss map according to Blaschke \cite{blaschke} (cf. \cite{bryant, rigoli}). 
\end{definition}

For a positive function $\lambda$ on the sphere $\mathbb{S}^3$ we consider the conformal metric $\lambda^2g_0$ on the sphere $\mathbb{S}^3$, which can be realized
as the 3-sphere $\mathbb{S}^3_\lambda$: $\lambda(1, x):\mathbb{S}^3\to\mathbb{N}^4_+\subset\mathbb{R}^{1,4}$ in Minkowski spacetime.  It is then very crucial and important to realize that the 
surface $\hat x$ in the 3-sphere $\mathbb{S}^3$ with the conformal metric $\lambda^2g_0$ is realized as the 2-surface $\hat\lambda (1, \hat x): \textup{M}^2\to \mathbb{N}^4_+
\subset \mathbb{R}^{1,4}$ inside the 3-sphere $\mathbb{S}^3_\lambda$. It is helpful to see the calculations in Appendix \ref{gauss} about the geometry of the 3-sphere 
$\mathbb{S}^3_\lambda$ in Minkowski spacetime $\mathbb{R}^{1,4}$.

\begin{lemma}\label{Lem:xi-lambda} If one works with a conformal metric $\lambda^2 g_0$ in general, then 
\begin{equation}\label{xi-lambda}
\xi = \xi_\lambda =  H_\lambda y_\lambda  + \overset{\rightarrow}{\bf n}_\lambda ,
\end{equation}
where $\overset{\rightarrow}{\bf n}_\lambda = \overset{\rightarrow}{\bf n} + (\log\lambda)_{\bf n} y$ is the unit normal to the surface 
$$y_\lambda = \hat \lambda (1, \hat x): \text{M}^2\to\mathbb{S}^3_\lambda\subset\mathbb{N}^4_+.$$
\end{lemma}

\proof It is easily seen that the normal direction to the surface $y_\lambda$ inside $\mathbb{S}^3_\lambda$ is $\lambda_{\bf n}(1, \hat x) + \lambda(0, {\bf n})$ and 
$< \lambda_{\bf n}(1, \hat x) + \lambda(0, {\bf n}), \lambda_{\bf n}(1, \hat x) + \lambda(0, {\bf n})> = \lambda^2$. Therefore the unit normal for the surface $y_\lambda$ in
$\mathbb{S}^3_\lambda$ is $\overset{\rightarrow}{\bf n}_\lambda = \overset{\rightarrow}{\bf n} + (\log\lambda)_{\bf n} y$. Then it is easily verified that
$$
H_\lambda y_\lambda + \overset{\rightarrow}{\bf n}_\lambda = Hy + \overset{\rightarrow}{\bf n}
$$
using \eqref{after-lambda}
\endproof

In the light of \eqref{mobius-metric},
the conformal Gauss map gives rise a spacelike 2-surface 
$$
\xi: \textup{M}^2\to \mathbb{S}^{1,3}\subset \mathbb{R}^{1,4}
$$
when the original surface $\hat x: \textup{M}^2\to\mathbb{S}^3$ is free of umbilical point. We will have more detailed discussions for the reasons to 
call $\xi$ the conformal Gauss map in Section \ref{canonicity}. \\

It is very interesting to see that Blaschke came across to the conformal Gauss map in a very different perspective. Blaschke considered the family of mean curvature 
2-spheres to the surface $\hat x$ in $\mathbb{S}^3$. A round 2-sphere in 3-sphere can be thought as the intersection of a timelike hyperplane and the 3-sphere
at time $t=1$ in Minkowski spacetime $\mathbb{R}^{1,4}$ and a timelike hyperplane in $\mathbb{R}^{1,4}$ is described by a unit normal vector lying in de Sitter 4-spacetime
$\mathbb{S}^{1,3}$. Given a direction $(H, H\hat x + {\bf n})\in \mathbb{S}^{1,3}$, the hyperplane perpendicular to that in $\mathbb{R}^{1,4}$ is given by the first equation 
in \eqref{enveloping}:
\begin{equation}\label{mean-curvature-sphere-1}
< (s,  z), (H, H\hat x + {\bf n})> = 0,
\end{equation} 
which is 
$$
-sH + H z\cdot (\hat x + \frac 1H {\bf n}) = 0.
$$
At the level $s=1$ in the 3-sphere $|z|=1$, we arrive at
$$
1 - \hat z \cdot (\hat x + \frac 1H {\bf n}) = 0.
$$
Then we may rewrite it as 
\begin{equation}\label{mean-curvature-sphere-2}
|\hat z - (\hat x + \frac 1H{\bf n})|^2 = \frac 1{H^2}
\end{equation}
which clearly is a round 2-sphere of mean curvature $H$ when intersects with the 3-sphere $\mathbb{S}^3\subset\mathbb{R}^4$ at $t=1$ in $\mathbb{R}^{1,4}$. Hence the 
equations \eqref{enveloping} exactly ask the surface
$y = (1, \hat x): \textup{M}^2\to \mathbb{S}^3\subset \mathbb{N}^4_+\subset \mathbb{R}^{1,4}$ is an envelope surface of the family of mean curvature 
2-spheres described by the conformal Gauss map $\xi$. \\

It is known that a mean curvature sphere of a surface goes to the mean curvature sphere of the image surface under conformal transformations.\\

\subsection{Conformal transforms}

Assume that the surface $\hat x: \textup{M}^2\to \mathbb{S}^3$ is free of umbilical point. Then the conformal Gauss map induces a spacelike 2-surface
in the de Sitter 4-space $\mathbb{S}^{1,3}$
$$
\xi: \textup{M}^2\to \mathbb{S}^{1,3}\subset \mathbb{R}^{1,4}.
$$
One notices that the equations \eqref{enveloping} imply that $y = (1, \hat x)$ is naturally a null normal vector the surface $\xi$ in the de Sitter 4-spacetime $\mathbb{S}^{1,3}$.
Because
$$
< y, \xi_{u^i}> = - <\xi, y_{u^i}> = 0.
$$
Hence it is natural to take the other null normal vector $y^*$ such that
\begin{equation}\label{y-star-def}
\aligned
<y^*, y> & = -1, \quad <y^*, y^*> = 0, \quad <y^*, \xi> = 0, \\ <y^*, \xi_{u^1}> & = 0, \text{ and } <y^*, \xi_{u^2}> = 0.\endaligned
\end{equation}
We may write 
$$
y^* = \hat\mu^* (1, \hat x^*).
$$

\begin{definition}\label{def-conf-trans} Suppose that $\hat x: \textup{M}^2\to\mathbb{S}^3$ is a surface with no umbilical point. And suppose that
$$
y^* =\hat \mu^*(1, \hat x^*): \textup{M}^2\to\mathbb{N}^4_+\subset \mathbb{R}^{1,4}
$$
satisfies the equations \eqref{y-star-def} for $y=(1, \hat x)$. Then the surface
$$
\hat x^*: \textup{M}^2\to \mathbb{S}^3
$$
is said to be the conformal transform of the surface $\hat x$ according to Robert Bryant \cite{bryant} (cf. \cite{blaschke}).
\end{definition}

It is important that the conformal transform $\hat x^*$ of a surface $\hat x$ is independent of the conformal factor $\lambda$. Notice that the equations in \eqref{y-star-def} remain
the same except the first one when replacing $y$ by $y_\lambda$.  It is again very interesting to recall how Blaschke 
discovered the surface $\hat x^*$. From the above discussions it is now easy to see that the surface $\hat x^*$ is nothing but the other envelope surface of the
family of round 2-spheres described by the conformal Gauss map $\xi$, i.e. the family of the mean curvature spheres of the surface $\hat x$. Since $y^*$ satisfies the last three
equations in \eqref{y-star-def}.\\

\subsection{The geometry of the surface $\xi$ in $\mathbb{S}^{1,3}$}\label{geometry-xi}
Recall that the first fundamental form for the surface $\xi$ in the de Sitter spacetime $\mathbb{S}^{1,3}\subset \mathbb{R}^{1,4}$ is
\begin{equation}\label{mobius-metric-1}
I^\xi = <d\xi, d\xi> = m|du|^2,
\end{equation}
where 
\begin{equation}\label{m-def}
m = \frac 12 E|\overset{\circ}{II}|^2.
\end{equation}
The first fundamental form $I^\xi$ is usually called the M\"{o}bius metric on the surface $\hat x$. We remark here that, if one works with a conformal metric $\lambda^2g_0$ instead,
then the M\"{o}bius metric remains the same
\begin{equation}\label{m-lambda}
m = m_\lambda = \frac 12 E_\lambda|\overset{\circ}{II}_\lambda|^2.
\end{equation}
The second fundamental form for the surface $\xi$ in $\mathbb{S}^{1,3}$ is 
given by
$$
II^{\xi} = - <d\xi, dy>  y - <d\xi, dy^*>  y^* = \Omega  y + \Omega^* y^* = \Omega_\lambda \hat\lambda^{-2} y_\lambda + \Omega^*_\lambda \hat\lambda^2 y^*_\lambda
$$
and 
\begin{equation}\label{omegas}
\aligned
\Omega_{ij}  = - <\xi_{u^i}, y_{u^j} >  &\text{ and } \ \Omega^*_{ij} = - <\xi_{u^i}, y^*_{u^j}>\\
(\Omega_\lambda)_{ij} = - <\xi_{u^i}, (y_\lambda)_{u^j} >  = \hat\lambda\Omega_{ij} & \text{ and } \ (\Omega^*_\lambda)_{ij} = - <\xi_{u^i}, (y^*_\lambda)_{u^j}> = \hat\lambda^{-1}
\Omega^*_{ij}\endaligned.
\end{equation}
In fact it is easy to calculate that
\begin{equation}\label{Omega}
\Omega =  \left[\begin{matrix} \frac {e-g}2 & f\\f & \frac {g-e}2\end{matrix}\right] = \overset{\circ}{II} 
\end{equation}

Let us first calculate the mean curvature in the $y^*$ direction. We notice that
$$ 
<\Delta_0\xi, \ y^*_\lambda> = ((\Omega^*_\lambda)_{11}+ (\Omega^*_\lambda)_{22})
$$
while
$$
<\Delta_0\xi, \ y_\lambda> = ((\Omega_\lambda)_{11} + (\Omega_\lambda)_{22}) = 0.
$$
Based on the calculations
$$
\aligned
<\Delta_0\xi, \  \xi> & =  -2m\\
<\Delta_0\xi,  \ \xi_{u^1}>  &  = \frac 12 m_{u^1} - \frac 12 m_{u^1} = 0\\
<\Delta_0\xi, \  \xi_{u^2}> & = - \frac 12 m_{u^2} + \frac 12 m_{u^2}  = 0.\endaligned
$$
we obtain
\begin{equation}\label{laplace-xi-1}
\Delta_0\xi = -((\Omega^*_\lambda)_{11}+ (\Omega^*_\lambda)_{22})y_\lambda -2m\xi=  (- ((\Omega^*_\lambda)_{11}+ (\Omega^*_\lambda)_{22}) -2mH_\lambda) y_\lambda
-2m\overset{\rightarrow}{\bf n}_\lambda  .
\end{equation}
On the other hand, we directly calculate
\begin{equation}\label{laplace-xi-2}
\aligned 
\Delta_0\xi  & = \Delta_0 (H_\lambda y_\lambda +\overset{\rightarrow}{\bf n}_\lambda) \\
& = (\Delta_0 H_\lambda) y_\lambda + H_\lambda \Delta_0 y_\lambda + 2 (H_\lambda)_{u^1} (y_\lambda)_{u^1} + 2 (H_\lambda)_{u^2} (y_\lambda)_{u^2} + 
\Delta_0\overset{\rightarrow}{\bf n}_\lambda \endaligned
\end{equation}  
It seems that the best way to calculate geometrically is to use the Lorentz orthogonal frame 
$$
\{y_\lambda,  y^\dagger_\lambda, (y_\lambda)_{u^1}, (y_\lambda)_{u^2}, \overset{\rightarrow}{\bf n}_\lambda\},
$$
where 
\begin{equation}\label{y-dagger-def}
\aligned
<y^\dagger_\lambda, y_\lambda> & = -1 \text{ and } \\
<y^\dagger_\lambda, y^\dagger_\lambda> & =  <y^\dagger_\lambda, (y_\lambda)_{u^1}> = <y^\dagger_\lambda, (y_\lambda)_{u^2}> 
= <y^\dagger_\lambda,  \overset{\rightarrow}{\bf n}_\lambda> = 0.\endaligned
\end{equation}
It is actually easy to find that 
\begin{equation}\label{y-dagger-explicit}
y_\lambda^\dagger = \frac 1{\lambda}(\frac 12|\nabla\log\lambda|^2 y + y^\dagger - \nabla \log\lambda),
\end{equation}
where $y^\dagger = \frac 12(1, - \hat x)$ and $\nabla$ is the gradient on the standard round 3-sphere. We will do inner product to both \eqref{laplace-xi-1} and \eqref{laplace-xi-2}
with the null vector $y^\dagger_\lambda$. To calculate 
$H_\lambda <\Delta_0 y_\lambda, y^\dagger_\lambda> + <\Delta_0\overset{\rightarrow}{\bf n}_\lambda, y^\dagger_\lambda>$ we rewrite
$$
H_\lambda <\Delta_0 y_\lambda, y^\dagger_\lambda>  = - H_\lambda (<(y_\lambda)_{u^1}, (y^\dagger_\lambda)_{u^1}> + <(y_\lambda)_{u^2}, (y^\dagger_\lambda)_{u^2}>)
$$
and
$$
<\Delta_0\overset{\rightarrow}{\bf n}_\lambda, y^\dagger_\lambda> = - <(\overset{\rightarrow}{\bf n}_\lambda)_{u^1}, (y^\dagger_\lambda)_{u^1}> - 
<(\overset{\rightarrow}{\bf n}_\lambda)_{u^2}, (y^\dagger_\lambda)_{u^2}>  - <\overset{\rightarrow}{\bf n}_\lambda, (y^\dagger_\lambda)_{u^i}>_{u^i} .
$$
Meanwhile one may calculate 
\begin{equation}\label{n-lambda-i-j}
\left\{\aligned
(\overset{\rightarrow}{\bf n}_\lambda)_{u^1} & = - \frac {e_\lambda}{E_\lambda} (y_\lambda)_{u^1} - \frac {f_\lambda}{E_\lambda} (y_\lambda)_{u^2} 
- <(\overset{\rightarrow}{\bf n}_\lambda)_{u^1}, y^\dagger_\lambda>y_\lambda \\
(\overset{\rightarrow}{\bf n}_\lambda)_{u^2} & = - \frac {f_\lambda}{E_\lambda} (y_\lambda)_{u^1} - \frac {g_\lambda}{E_\lambda} (y_\lambda)_{u^2}
- <(\overset{\rightarrow}{\bf n}_\lambda)_{u^2}, y^\dagger_\lambda>y_\lambda.
\endaligned\right.
\end{equation}
Hence we have
\begin{equation}\label{trace-Omega-star}
\aligned
H_\lambda &  <\Delta_0 y_\lambda, y^\dagger_\lambda> + <\Delta_0\overset{\rightarrow}{\bf n}_\lambda, y^\dagger_\lambda>\\ &  = 
E^{-1}_\lambda (\overset{\circ}{II}_\lambda)_{ij} <(y_\lambda)_{u^i}, (y^\dagger_\lambda)_{u^j}>  - <\overset{\rightarrow}{\bf n}_\lambda, (y^\dagger_\lambda)_{u^i}>_{u^i} \\
& = - E_\lambda^{-1} (\overset{\circ}{II}_\lambda)_{ij}
R^\lambda_{i3j3} + E_\lambda (R^\lambda)_{3i,}^{\quad i}
\endaligned
\end{equation}
due to \eqref{n-dagger}, \eqref{i-dagger-j}, and \eqref{coord-covar}.  Now we obtain the mean curvature of the surface $\xi$ in the de Sitter spacetime $\mathbb{S}^{1,3}$.

\begin{lemma} Suppose that $\hat x:\textup{M}^2\to\mathbb{S}^3$ is an immersed surface with no umbilical point 
and that $\xi:\textup{M}^2\to \mathbb{S}^{1,3}$ is the conformal Gauss map. Then the surface $\xi$ is spacelike and its mean curvature is a null vector
\begin{equation}\label{mean-curvature-xi}
H^\xi  = 2\hat\lambda^2\frac {{\mathcal H}_\lambda}{|\overset{\circ}{II}_\lambda|^2} y^*_\lambda 
\end{equation}
for any positive function $\lambda$ on the 3-sphere $\mathbb{S}^3$, where 
\begin{equation}\label{script-H-lambda}
{\mathcal H}_\lambda = \Delta_\lambda H_\lambda  + |\overset{\circ}{II}_\lambda|^2H_\lambda
+ (\overset{\circ}{II}_\lambda)^{ij}(R^\lambda)_{i3j3}  -  (R^\lambda)_{3i,}^{\ \ i}, 
\end{equation}
$(R^\lambda)_{ijkl}$ and $(R^\lambda)_{ij}$  are the Riemann curvature and Ricci curvature for the conformal metric $\lambda^2g_0$ on the 3-sphere
$\mathbb{S}^3$ respectively.
\end{lemma}
\proof We perform inner product  to \eqref{laplace-xi-1} and \eqref{laplace-xi-2}  by the null vector $y^\dagger_\lambda$ and obtain that
\begin{equation}\label{trace-star}
(\Omega^*_\lambda)_{11} + (\Omega^*_\lambda)_{22}  = E_\lambda (- \Delta_\lambda H_\lambda - |\overset{\circ}{II}_\lambda|^2H_\lambda - 
(\overset{\circ}{II}_\lambda)^{ij}(R^\lambda)_{i3j3}  + (R^\lambda)_{3i,}^{\ \ i})
\end{equation}
in the light of  \eqref{trace-Omega-star}. Then one can easily calculate the mean curvature for $\xi$ in $\mathbb{S}^{1,3}$.
\endproof 

\vskip 0.2in
We remark that \eqref{mean-curvature-xi} actually shows that 
\begin{equation}\label{2-scalar-invar}
{\mathcal H}_\lambda = \hat\lambda^{-3} (- \Delta H - |\overset{\circ}{II}|^2H)
\end{equation}
for a surface $\hat x$ in the conformal 3-sphere.

\begin{theorem} \label{conf-gauss-surface} (\cite{blaschke} \cite{bryant}) Suppose that $\hat x: \textup{M}^2\to\mathbb{S}^3$ is an immersed surface with no
umbilical point. Then $\hat x$ is a Willmore surface in $\mathbb{S}^3$ if and only if the conformal Gauss map induces a minimal spacelike surface in 
the de Sitter spacetime $\mathbb{S}^{1,3}$. Moreover its conformal transform $\hat x^*$ is a dual Willmore surface in $\mathbb{S}^3$.
\end{theorem}

\proof  Most of this theorem has been known to Blaschke \cite{blaschke} and Bryant \cite{bryant}. Because 
Lemma \ref{xi-2} implies that $\xi$ is also the conformal Gauss map (up to the sign) for $\hat x^*$ when $H^\xi$ vanishes. 
The two dual Willmore surfaces are the two envelope surfaces of the family of round 2-spheres described by the conformal Gauss map $\xi$.
\endproof

\begin{remark} It is also known to Balschke \cite{blaschke} and Bryant \cite{bryant} that
\begin{itemize}
\item If $\hat x$ is a minimal surface in $\mathbb{S}^3$, then $\hat x^* = -\hat x$.
\item $\hat x$ is a Willmore surface if and only if $\hat x^{**} = \hat x$,
which raises an interesting question: what does it mean $\hat x^{***} = \hat x$ if possible?
\end{itemize}

\end{remark} 

\subsection{Finding $y^*_\lambda$}
Let us now solve $y^*_\lambda$ for $y_\lambda = \hat\lambda (1, \hat x) = \hat\lambda y$, where $\hat\lambda = \lambda\circ \hat x$ and $\lambda$ is 
a positive function  on the sphere $\mathbb{S}^3$. At each point on the surface we set
$$
y^*_\lambda =\kappa y_\lambda +  \kappa_\dagger y^\dagger_\lambda + b \overset{\rightarrow}{\bf n}_\lambda + \frac {\omega^\lambda_1}{E_\lambda} (y_\lambda)_{u^1} 
+\frac {\omega^\lambda_2}{E_\lambda} (y_\lambda)_{u^2}.
$$
And we get from \eqref{y-star-def}
\begin{equation}\label{abcd}
\left\{ \aligned \kappa_\dagger & = 1\\
- 2\kappa_+\kappa_- + b^2 + \frac {(\omega^\lambda_1)^2 + (\omega^\lambda_2)^2}{E_\lambda} & = 0 \\
b & = H_\lambda\\
- (\Omega_\lambda)_{11}\omega^\lambda_1 - (\Omega_\lambda)_{12} \omega^\lambda_2 & = (H_\lambda)_{u^1}E_\lambda\\
- (\Omega_\lambda)_{21}\omega^\lambda_1 - (\Omega_\lambda)_{22} \omega^\lambda_2 & = (H_\lambda)_{u^2} E_\lambda.\endaligned\right.
\end{equation}
We therefore have

\begin{lemma}\label{lemma-x-star} Suppose that $\hat x: \textup{M}^2\to\mathbb{S}^3$ is an immersed surface with no umbilical point.
Then
\begin{equation}\label{y-star-lambda}
 y^*_\lambda  =\frac 12  ( |\omega^\lambda|^2 + H_\lambda^2)y_\lambda + y_\lambda^\dagger 
 + H_\lambda\overset{\rightarrow}{\bf n}_\lambda  - (\overset{\circ}{II})^{-1}_\lambda d H_\lambda
\end{equation}
for any positive function $\lambda$ on the 3-sphere, where 
$$
|\omega^\lambda|^2 =  \frac {(\omega^\lambda_1)^2 + (\omega^\lambda_2)^2}{E_\lambda}  = \frac 1m((H_\lambda)_{u^1}^2 +(H_\lambda)_{u^2}^2).  
$$
In particular,
\begin{equation}\label{y-star}
y^*  = \frac 12 (|\omega|^2 +H^2) y + \frac 12(1, -\hat x)  + H (0, {\bf n})   - (0,  (\overset{\circ}{II})^{-1} d H),
\end{equation}
and
\begin{equation}\label{x-star}
x^*  = a \hat x + \frac H {1-a} {\bf n}   - \frac 1{1-a}(\overset{\circ}{II})^{-1} d H,
\end{equation}
where 
\begin{equation}\label{what-a}
a = \frac {|\omega|^2 +H^2-1}{|\omega|^2 +H^2 +1}.
\end{equation}
\end{lemma}
\proof One simply solves \eqref{abcd} if $\det\Omega_\lambda \neq 0$, which is equivalent to the fact that the surface has no umbilical point.
\endproof

\subsection{Canonicity of $y^*$}\label{canonicity}

Now we want to show that the choice of $y^*$ is canonical in terms of doing conformal geometry for the surface $\hat x$ in $\mathbb{S}^3$. 
It is important to realize that there are two separate issues here. One is about the symmetry of the conformal 3-sphere. To be precise, for a conformal transformation
$$
\phi: \mathbb{S}^3\to\mathbb{S}^3
$$
and the transformed surface $$\phi(\hat x): \textup{M}^2\to \mathbb{S}^3,$$ is it true that 
$$
\tilde\phi (\tilde x)= \alpha\tilde\phi (y) + \alpha\rho\tilde\phi(y^*):\mathbb{R}^+\times\mathbb{R}^+\times\textup{M}^2\to\mathbb{R}^{1,4}$$ 
is the associate 4-surface of $\phi(\hat x)$  in $\mathbb{R}^{1,4}$, where $\tilde\phi$ is the corresponding Lorentz transformation
on $\mathbb{R}^{1,4}$ to $\phi$? 
The other issue is whether or not the associate surface $\tilde x$ is independent of metrics in the conformal class of the round 3-sphere.
The first easy and important fact is that the conformal Gauss map is independent of the metrics in the conformal class.

\begin{lemma}\label{transform-gauss} Suppose that $\hat x: \textup{M}^2\to\mathbb{S}^3$ is an immersed surface. Then the conformal Gauss map $\xi$ is independent of the metrics
in the conformal class of the round 3-sphere $\mathbb{S}^3$. 
Meanwhile, the conformal Gauss map for the transformed surface $\phi(\hat x)$
 is exactly $\tilde\phi(\xi)$, where $\tilde\phi$ is the Lorentz transformation on the Minkowski spacetime
$\mathbb{R}^{1,4}$ corresponding to a conformal transformation $\phi$ on $\mathbb{S}^3$. 
\end{lemma}

\proof First of all, one needs to realize that, for any given metric in the conformal class of the round 3-sphere, it simply amounts to consider the surface
$$
y_\lambda = \hat \lambda(1, \hat x): \textup{M}^2\to \mathbb{N}^4_+
$$
for some positive function $\lambda: \mathbb{S}^3\to \mathbb{R}^+$ and $\hat\lambda = \lambda\circ\hat x$. 
But this only possibly alters the parametrization of the homogeneous null 3-surface 
$x^{\mathbb{N}} = \alpha \hat \lambda (1, \hat x): \mathbb{R}^+\times\textup{M}^2\to \mathbb{N}^4_+$. Hence it will not alter the conformal Gauss map. 
Of course one has already seen this from Lemma \ref{Lem:xi-lambda}.\\

Next we consider the transformed surface $\phi(\hat x)$. Recall that, given a conformal transformation $\phi$ of 3-sphere, we have a unique Lorentz transformation
$\tilde\phi$ in the time and orientation preserving component of the Lorentz group on the Minkowski spacetime such that, for $\lambda(1, \hat x)\in \mathbb{R}^{1, 4}$, 
\begin{equation}\label{poincare}
\tilde \phi (\lambda (1, \hat x)) = \lambda \mu (1, \phi(\hat x))
\end{equation} 
for some positive number $\mu$. By the definition, which requires $\tilde\phi$ is a linear map and
$$
<\tilde\phi ((t, \hat x)), \tilde\phi((s, \hat y))> = < (t, \hat x), (s, \hat y)>,
$$
we now easily see that $\tilde\phi(\xi)$ is the conformal Gauss map for the transformed surface $\phi(\hat x)$. Since $\tilde\phi(\xi)$ is the unit normal vector field to the 
homogeneous null 3-surface
$\tilde\phi(x)$ in $\mathbb{N}^4_+$ that is conformal map from $\textup{M}^2$ to $\mathbb{S}^{1,3}$. 
\endproof

Consequently we have 

\begin{proposition} Suppose that $\hat x: \textup{M}^2\to\mathbb{S}^3$ is an immersed surface with no umbilical point. Then the associate surface
$$
\tilde x = \alpha y_\lambda + \alpha\rho y^*_\lambda: \mathbb{R}^+\times\mathbb{R}^+\times\textup{M}^2\to\mathbb{R}^{1,4},
$$
for any $y_\lambda= \hat\lambda (1, \hat x)$ and $y^* = \hat\lambda^{-1} \lambda^*(1, \hat x^*)$ defined by the equations \eqref{y-star-def}, is independent of the metrics in 
conformal class of the round 3-sphere $\mathbb{S}^3$.
\end{proposition}

\proof It suffices to verify that 
\begin{equation}\label{star-scaled}
(\hat\lambda y)^* = \hat\lambda^{-1} y^*.
\end{equation}
Since it implies that the change of metrics in the conformal class will at most cause possible change of parametrization of the associate surface $\tilde x$. 
\endproof

We also have from Lemma \ref {transform-gauss} the following:

\begin{lemma} Suppose that $\hat x: \textup{M}^2\to\mathbb{S}^3$ is an immersed surface with no umbilical point. Let 
$y_\lambda = \hat\lambda (1, \hat x)\in \mathbb{N}^4_+$ and let $\phi$ be a conformal transformation of 3-sphere. Then 
\begin{equation}\label{transform-y-star}
\tilde\phi(y_\lambda)^* = \tilde\phi(y^*_\lambda).
\end{equation}
Hence 
\begin{equation} \label{transform-x-star} 
\phi(\hat x^*) = (\phi(\hat x))^*.
\end{equation}
\end{lemma}
\proof From Lemma \ref{transform-gauss} we know that the conformal Gauss map for the transformed surface $\phi(\hat x)$ is $\tilde\phi(\xi)$. Then it is easy to verify 
\eqref{y-star-def} for $\tilde\phi(y^*)$ to be $\tilde\phi(y)^*$. Then the equation \eqref{transform-x-star} follows from \eqref{poincare} and \eqref{transform-y-star}:
$$
\hat\gamma^*(1, (\phi(\hat x))^*) = \tilde\phi(y)^* = \tilde\phi(y^*) = \hat\mu^*\hat\lambda^*(1, \phi(\hat x^*)).
$$
\endproof

Therefore we have 

\begin{proposition} Suppose that $\hat x: \textup{M}^2\to\mathbb{S}^3$ is an immersed surface with no umbilical point. Let 
$\phi$ be a conformal transformation of 3-sphere. Then the associate 4-surface in $\mathbb{R}^{1,4}$ of the 
transformed surface $\phi(\hat x)$ is exactly the 4-surface $\tilde\phi(\tilde x)$ transformed from the associate 4-surface $\tilde x$ 
of the original surface $\hat x$ under the corresponding Lorentz transformation $\tilde\phi$ of $\phi$.
\end{proposition}

\section{The geometry of the associate surfaces}

In this section we calculate the first and second fundamental forms for the associate homogeneous timelike 4-surfaces $\tilde x$ in 
$\mathbb{R}^{1,4}$ as well as for the associate ruled surface $x^+$ in the hyperboloid $\mathbb{H}^4$,
for a given immersed 2-surface $\hat x$  in $\mathbb{S}^3$.

\subsection{The first fundamental form for $\tilde x$ in $\mathbb{R}^{1,4}$}\label{Sect:I-tilde}

To calculate the first fundamental form for the surface in the parametrization 
\begin{equation}\label{lambda-para}
\tilde x = \alpha y_\lambda + \alpha \rho y^*_\lambda 
\end{equation}
associated with a conformal metric $\lambda^2g_0$ on the 3-sphere $\mathbb{S}^3$ , we first calculate
$$
d\tilde x  = (y_\lambda +\rho y^*_\lambda)d\alpha + \alpha y^*_\lambda d\rho + (\alpha (y_\lambda)_{u^1}+\alpha\rho (y^*_\lambda)_{u^1})du^1 
+ (\alpha (y_\lambda)_{u^2} + \alpha\rho (y^*_\lambda)_{u^2})du^2.
$$
Hence the first fundamental form for the associate 4-surface $\tilde x$ in the coordinates $(\alpha, \rho, u^1, u^2)$ is 
$$
\aligned
I^{\tilde x} & = <d\tilde x, d\tilde x> = -2 \rho d\alpha d\alpha - 2\alpha d\alpha d\rho \\ 
& + 2\alpha^2 <(y^*_\lambda, (y_\lambda)_{u^1}>d\rho du^1 + 2\alpha^2 <y^*_\lambda, (y_\lambda)_{u^2}>d\rho du^2\\
& + <\alpha (y_\lambda)_{u^1}+\alpha\rho (y^*_\lambda)_{u^1}, \alpha (y_\lambda)_{u^1}+\alpha\rho (y^*_\lambda)_{u^1}>(du^1)^2 \\
& +  <\alpha (y_\lambda)_{u^2} +\alpha\rho (y^*_\lambda)_{u^2}, \alpha (y_\lambda)_{u^2}+\alpha\rho (y^*_\lambda)_{u^2}>(du^2)^2\\
& + 2<\alpha (y_\lambda)_{u^1}+\alpha\rho (y^*_\lambda)_{u^1}, \alpha (y_\lambda)_{u^2}+\alpha\rho (y^*_\lambda)_{u^2}>du^1du^2.
\endaligned
$$
In fact one may calculate 
\begin{equation}\label{tilde-x-I}
\left\{ \aligned (y_\lambda)_{u^1} & =  - \omega^\lambda_1 y_\lambda - \frac {(\Omega_\lambda)_{11} }{m}\xi_{u^1} - \frac {(\Omega_\lambda)_{12}}{m}\xi_{u^2}\\
(y_\lambda)_{u^2} & =  - \omega^\lambda_2 y_\lambda - \frac {(\Omega_\lambda)_{21} }{m}\xi_{u^1} - \frac {(\Omega_\lambda)_{22}}{m}\xi_{u^2}\\
(y^*_\lambda)_{u^1} & =   \omega^\lambda_1 y^*_\lambda - \frac {(\Omega_\lambda)^*_{11}}{m}\xi_{u^1} - \frac {(\Omega_\lambda)^*_{12}}{m}\xi_{u^2}\\
(y^*_\lambda)_{u^2} & =   \omega^\lambda_2 y^*_\lambda - \frac {(\Omega_\lambda)^*_{21} }{m}\xi_{u^1} - \frac {(\Omega_\lambda)^*_{22}}{m}\xi_{u^2}\endaligned
\right.
\end{equation}
where
\begin{equation}\label{omega}
\omega^\lambda = <dy_\lambda, y^*_\lambda>  = - I_\lambda(\Omega_\lambda^{-1}dH_\lambda)
\end{equation}
based on \eqref{abcd}. Now let us write $I^{\tilde x}$ in matrix form:
\begin{equation}\label{matrix-tilde}
I_{\tilde x} = \left[\begin{matrix} \begin{matrix} \ -2\rho & -\alpha \\ -\alpha & \ 0\end{matrix} & \begin{matrix} 0 & 0 \\
\alpha^2\omega^\lambda_1 & \alpha^2 \omega^\lambda_2\end{matrix}\\
\begin{matrix} \quad\quad 0 & \quad \alpha^2\omega^\lambda_1\\ \quad\quad 0 & \quad \alpha^2\omega^\lambda_2\end{matrix} & \quad 
\alpha^2 F \quad\end{matrix}\right]
 \end{equation}
 where
 \begin{equation}\label{F-matrix}
 \left\{
 \aligned
 F_{11} & = \frac 1m(p^2 + q^2)  +2 \rho (\omega^\lambda_1)^2\\
 F_{12} & =  F_{21} =  \frac 1m q(p+r)  + 2 \rho \omega^\lambda_1\omega^\lambda_2\\
  F_{22} & = \frac 1m(q^2 + r^2)  + 2 \rho(\omega^\lambda_2)^2 \endaligned\right.
\end{equation}
 and
$$
\left[\begin{matrix} p & q \\ q & r\end{matrix}\right] = \Omega_\lambda + \rho\Omega^*_\lambda.
$$
It can be calculated that 
\begin{equation}\label{det-I-tilde}
\det I^{\tilde x} = - \frac {\alpha^6}{m^2}(pr - q^2)^2 = -\frac {\alpha^6}{4m^2} (E^2_\lambda|\Omega_\lambda + \rho\Omega^*_\lambda|^2 
- \rho^2((\Omega^*_\lambda)_{11}+ (\Omega^*_\lambda)_{22})^2)^2
\end{equation}
which can tell us where the associate surface $\tilde x$ is degenerate. It is maybe a little surprising that it is actually not difficult to calculate the inverse of 
$I_{\tilde x}$. We present the calculations in Appendix \ref{app-inverse-tilde} since they are straightforward calculations.
  
 \subsection{The second fundamental form for $\tilde x$ in $\mathbb{R}^{1,4}$}\label{Sect:II-tilde}
 
 It is clear from the definition that the conformal Gauss map $\xi$ is the unit normal vector for the associate 4-surface $\tilde x$ in $\mathbb{R}^{1,4}$. Hence 
 the second fundamental form for $\tilde x$ in $\mathbb{R}^{1,4}$ is 
 \begin{equation}\label{tilde-x-II}
 II^{\tilde x} = - < d\tilde x, d\xi> = (\alpha(\Omega_\lambda)_{ij}+\alpha\rho(\Omega^*_\lambda)_{ij})du^idu^j
 \end{equation}
 or in matrix form
 $$
II_{\tilde x} =  \left[\begin{matrix} \quad 0 & 0 \\ \quad 0 & \alpha\Omega_\lambda +\alpha \rho\Omega^*_\lambda\end{matrix}\right].
$$
Therefore the mean curvature for the associate 4-surface in $\mathbb{R}^{1,4}$ is
$$
 H^{\tilde x} = \text{Tr} (I_{\tilde x})^{-1}II_{\tilde x}.
 $$
 To calculate the mean curvature $H^{\tilde x}$ one only needs to know the low-right $2\times 2$ block 
 in the inverse of the matrix $I_{\tilde x}$. According to the calculations in Appendix \ref{app-inverse-tilde},  particularly \eqref{F-star-inverse} \eqref{app-3j} \eqref{app-4j}, 
 we therefore have
\begin{equation}\label{x-tilde-mean}
\aligned
H^{\tilde x} & = \frac m{\alpha(pr-q^2)^2}((q^2+r^2)p - 2 q^2(p+r) + (p^2+q^2)r)\\
& = \frac {m(p+r)}{\alpha(pr-q^2)},  \endaligned
\end{equation}
where
$$
pr - q^2 =  \det\Omega_\lambda - \rho\text{Tr}\Omega_\lambda\Omega^*_\lambda + \rho^2\det\Omega^*_\lambda
$$
and
\begin{equation}\label{script-H-0}
p+r  =  \rho((\Omega^*_\lambda)_{11}+(\Omega^*_\lambda)_{22})   =  - \rho  E_\lambda {\mathcal H}_\lambda
\end{equation}
in the light of \eqref{trace-star}.

\begin{theorem}\label{Th:tilde-H} Suppose that $\hat x:\textup{M}^2\to\mathbb{S}^3$ is an immersed surface with no umbilical point. Then $\hat x$ is a Willmore surface
in $\mathbb{S}^3$ if and only if the associate 4-surface $\tilde x$ in $\mathbb{R}^{1,4}$ is minimal. 
\end{theorem}
\proof Based on the above equations \eqref{script-H-0} and \eqref{x-tilde-mean} we obtain that
\begin{equation}\label{script-H}
H^{\tilde x} = \frac {\rho \det\Omega_\lambda {\mathcal H}_\lambda}
{\alpha (\det\Omega_\lambda - \rho\text{Tr}\Omega_\lambda\Omega^*_\lambda + \rho^2\det\Omega^*_\lambda)}.
\end{equation}
\endproof

\subsection{Local fundamental theorem for surfaces in conformal 3-sphere} 

In this subsection we want to state and prove a local fundamental theorem for surfaces in conformal 3-sphere. In the previous section we have introduced
the associate surface $\tilde x$ in Minkowski spacetime $\mathbb{R}^{1, 4}$ from a given surface $\hat x$ in $\mathbb{S}^3$. From the geometric structure of the associate surface
$\tilde x$ one can tell that its intersection with the positive light cone $\mathbb{N}^4_+$ is a homogeneous null 3-surface whose projectivization will recover the original surface 
$\hat x$ in $\mathbb{S}^3$. 

Given a surface $\hat x$ in $\mathbb{S}^3$ with a isothermal coordinates $(u^1, u^2)$ on the parameter space $\textup{M}^2$, we have the first fundamental form $I$ in matrix form
$$
I = \left[\begin{matrix} E & 0\\0 & E\end{matrix}\right]
$$
and
the second fundamental fundamental $II$ form in matrix form
$$
II = \left[\begin{matrix} e & f\\f & g\end{matrix}\right]
$$
The local fundamental theorem for surfaces in Riemannian geometry states that, up to isometries of the standard round sphere $\mathbb{S}^3$, 
locally the surface is uniquely determined by the first fundamental form $I$ and the second fundamental form $II$ in the standard round sphere $\mathbb{S}^3$. 
Conversely, given a positive definite symmetric 2-form $I$ and a symmetric 2-form $II$ in the parameter domain, which satisfy some integrability conditions 
(Gauss-Codazzi equations), up to isometries,  there is locally a unique surface $\hat x$ in the standard round sphere $\mathbb{S}^3$ whose first and 
second fundamental forms are $I$ and $II$. We are looking for the analogous local fundamental theorem for surfaces in conformal round 3-sphere $\mathbb{S}^3$. 
The core idea of the local fundamental theorem in Riemannian geometry is to solve the structure equations, which are the equations of motion of Fren\'{e}t frames on the
surface and are determined from $I$ and $II$.  \\

Our strategy here is to use the local fundamental theorem for the associate surface $\tilde x$ in the Minkowski spacetime $\mathbb{R}^{1, 4}$ to establish the local fundamental 
theorem for 
a surface $\hat x$ in the conformal sphere $\mathbb{S}^3$. Since the association introduced in previous subsections requires that the surface $\hat x$ has no umbilical
point, we will always assume here that surfaces $\hat x$ have no umbilical point. \\

To summarize the previous discussions, given a surface $\hat x$ in $\mathbb{S}^3$, 
we have $I=E|du|^2$ and $II = e (du^1)^2 + 2f du^1du^2 + g(du^2)^2$. We also have the so-called M\"{o}bius metric $I^\xi = m |du|^2 = \frac 12 E|\overset{\circ}{II}|^2|du|^2$ 
induced from the Conformal Gauss map $\xi$ of the surface $\hat x$, where 
$$
\overset{\circ}{II}= \left[\begin{matrix} \frac {e-g}2  & f\\ f & \frac {g-e}2\end{matrix}\right]
$$ 
is the traceless part of the second fundamental form $II$. We then construct the associate surface
$$
\tilde x = \alpha y_\lambda +\alpha \rho y^*_\lambda: \mathbb{R}^+\times\mathbb{R}^+\times\textup{M}^2: \mathbb{R}^{1, 4}.
$$
\\

The first fundamental form $I^{\tilde x}$ for $\tilde x$ in matrix form is, from 
\eqref{matrix-tilde},
$$
\left[\begin{matrix} \begin{matrix} -2\rho & -\alpha\\-\alpha & 0 \end{matrix} & \begin{matrix} 0 & 0 \\
\alpha^2\omega^\lambda_1 & \alpha^2\omega^\lambda_2\end{matrix} \\ \begin{matrix} \quad\quad 0 & \alpha^2\omega^\lambda_1\\
\quad\quad 0 & \alpha^2\omega^\lambda_2\end{matrix} & \begin{matrix} \frac {\alpha^2}m(p^2+q^2) + 2\alpha^2\rho(\omega^\lambda_1)^2 
& \frac {\alpha^2}m q(p+r)+ 2\alpha^2\rho\omega^\lambda_1\omega^\lambda_2\\ \frac {\alpha^2}m q(p+r) +2\alpha^2\rho\omega^\lambda_1\omega^\lambda_2& 
\frac {\alpha^2}m(q^2 +r^2) +2\alpha^2\rho(\omega^\lambda_2)^2\end{matrix} \end{matrix}\right],
$$
where the 1-form
$$
\omega^\lambda  = \omega^\lambda_1du^1 + \omega^\lambda_2du^2 = -d\log\hat\lambda - I(\Omega^{-1} (d H)) = d\log\hat\lambda +\omega.
$$
\\

And the second fundamental form $II^{\tilde x}$ for $\tilde x$ in $\mathbb{R}^{1, 4}$ in matrix form is, from \eqref{tilde-x-II},
$$
\left[\begin{matrix} 0 & 0 \\ 0 & \alpha \Omega_\lambda + \alpha\rho\Omega^*_\lambda
\end{matrix}\right].
$$
where $\Omega_\lambda = \hat\lambda\Omega$ and $\Omega^*_\lambda = \hat\lambda^{-1}\Omega^*$.
Notice that $I^{\tilde x}$ and $II^{\tilde x}$ are exactly determined by the M\"{o}bius metric $I^\xi = m|du|^2$, the 1-form $\omega$,  
the traceless symmetric 2-tensor $\Omega$ and the symmetric 2-tensor $\Omega^*$, plus the conformal factor $\hat \lambda$. \\

Next we write the equations for the motion of the Fren\'{e}t frames on the associate surface $\tilde x$ according to $I^{\tilde x}$ and $II^{\tilde x}$. We consider the Fren\'{e}t frame 
$\{y_\lambda, y^*_\lambda, \frac 1{\sqrt m}\xi_{u^1}, \frac 1{\sqrt m}\xi_{u^2}, \xi\}$ on the associate surface $\tilde x$. 
Because they are the orthonormal frames on $\tilde x$ with respect to the 
Minkowski metric $\tilde{\mathcal G}_0$ on $\mathbb{R}^{1, 4}$. We now write
\begin{equation}\label{structure-equation-1}
\frac \partial{\partial u^1}\left[\begin{matrix} y_\lambda\\y^*_\lambda\\ \frac 1{\sqrt m}\xi_{u^1}\\ \frac 1{\sqrt m}\xi_{u^2}\\ \xi\end{matrix}\right] = 
\left[\begin{matrix} -\omega^\lambda_1 & 0 & - \frac 1{\sqrt m}(\Omega_\lambda)_{11} & -\frac 1{\sqrt m}(\Omega_\lambda)_{12} & 0\\
0 & \omega^\lambda_1 & - \frac 1m(\Omega^*_\lambda)_{11} & - \frac 1m (\Omega^*_\lambda)_{12} & 0\\ 
\frac 1{\sqrt m}(\Omega_\lambda)_{11} &  \frac 1{\sqrt m}(\Omega^*_\lambda)_{11} & 0 &  - \frac 1{2m} m_{u^2} & -\sqrt{m}\\
\frac 1{\sqrt m}(\Omega_\lambda)_{21} &  \frac 1{\sqrt m}(\Omega^*_\lambda)_{21} & \frac 1{2m} m_{u^2} & 0 & 0 \\ 
0 & 0 & \sqrt{m} & 0 & 0\end{matrix}\right] \left[\begin{matrix} y_\lambda \\y^*_\lambda \\ \frac 1{\sqrt m}\xi_{u^1}\\ \frac 1{\sqrt m}\xi_{u^2}\\ \xi\end{matrix}\right] 
\end{equation}
and 
\begin{equation}\label{structure-equation-2}
\frac \partial{\partial u^2}\left[\begin{matrix} y_\lambda\\y^*_\lambda\\ \frac 1{\sqrt m}\xi_{u^1}\\ \frac 1{\sqrt m}\xi_{u^2}\\ \xi\end{matrix}\right] = 
\left[\begin{matrix} -\omega^\lambda_2 & 0 & - \frac 1{\sqrt m}(\Omega_\lambda)_{21} & -\frac 1{\sqrt m}(\Omega_\lambda)_{22} & 0\\
0 & \omega^\lambda_2 & - \frac 1m(\Omega^*_\lambda)_{21} & - \frac 1m (\Omega^*_\lambda)_{22} & 0\\ 
\frac 1{\sqrt m}(\Omega_\lambda)_{21} &  \frac 1{\sqrt m}(\Omega^*_\lambda)_{21} & 0 &  - \frac 1{2m} m_{u^1} & 0\\
\frac 1{\sqrt m}(\Omega_\lambda)_{22} &  \frac 1{\sqrt m}(\Omega^*_\lambda)_{22} & \frac 1{2m} m_{u^1} & 0 & -\sqrt{m} \\ 
0 & 0 & 0 &  \sqrt{m} & 0\end{matrix}\right] \left[\begin{matrix} y_\lambda \\y^*_\lambda \\ \frac 1{\sqrt m}\xi_{u^1}\\ \frac 1{\sqrt m}\xi_{u^2}\\ \xi\end{matrix}\right] 
\end{equation}
Remember we also have the two trivial equations 
$$
\frac \partial{\partial \alpha}\left[\begin{matrix} y_\lambda\\y^*_\lambda\\ \frac 1{\sqrt m}\xi_{u^1}\\ \frac 1{\sqrt m}\xi_{u^2}\\ \xi\end{matrix}\right]= 0 \text{ and } \frac \partial{\partial \rho}\left[\begin{matrix} y_\lambda\\y^*_\lambda\\ \frac 1{\sqrt m}\xi_{u^1}\\ \frac 1{\sqrt m}\xi_{u^2}\\ \xi\end{matrix}\right]= 0.
$$
To solve the systems \eqref{structure-equation-1} and \eqref{structure-equation-2} of ODE, the necessary integrable condition is
\begin{equation}\label{integrable}
\frac \partial{\partial u^1}\frac \partial{\partial u^2}\left[\begin{matrix} y_\lambda\\y^*_\lambda\\ \frac 1{\sqrt m}\xi_{u^1}\\ \frac 1{\sqrt m}\xi_{u^2}\\ \xi\end{matrix}\right]
= \frac \partial{\partial u^2}\frac\partial{\partial u^1}\left[\begin{matrix} y_\lambda\\y^*_\lambda\\ \frac 1{\sqrt m}\xi_{u^1}\\ \frac 1{\sqrt m}\xi_{u^2}\\ \xi\end{matrix}\right].
\end{equation}
It turns out \eqref{integrable} is equivalent to the following six equations on the variables: the positive function $m$, the 1-form $\omega^\lambda$, the traceless symmetric 
matrix $\Omega_\lambda$ and the symmetric matrix $\Omega^*_\lambda$, 
\begin{equation}\label{codazzi-y}
\left\{\aligned (\Omega_\lambda)_{11, 2} - (\Omega_\lambda)_{12,1} & = \omega^\lambda_1(\Omega_\lambda)_{12} - \omega^\lambda_2(\Omega_\lambda)_{11}\\
 (\Omega_\lambda)_{12,2} - (\Omega_\lambda)_{22,1} & = \omega^\lambda_1(\Omega_\lambda)_{22} - \omega^\lambda_2(\Omega_\lambda)_{12}\endaligned\right.
 \end{equation}
\begin{equation}\label{codazzi-y-star}
\left\{\aligned (\Omega^*_\lambda)_{11, 2}- (\Omega^*_\lambda)_{12,1}& = -  \omega^\lambda_1(\Omega^*_\lambda)_{12} + \omega^\lambda_2(\Omega^*_\lambda)_{11} 
+ \frac 12 \frac{(\Omega^*_\lambda)_{11} + (\Omega^*_\lambda)_{22}}{|\Omega_\lambda|^2} (|\Omega_\lambda|^2)_{u^2}\\
 (\Omega^*_\lambda)_{12,2}- (\Omega^*_\lambda)_{22,1}& = - \omega^\lambda_1(\Omega^*_\lambda)_{22} + \omega^\lambda_2(\Omega^*_\lambda)_{12}  
 + \frac 12 \frac{(\Omega^*_\lambda)_{11} + (\Omega^*_\lambda)_{22}}{|\Omega_\lambda|^2} (|\Omega_\lambda|^2)_{u^2}\endaligned\right.
 \end{equation}
 \begin{equation}\label{codazzi-mix}
 \omega^\lambda_{1, 2} - \omega^\lambda_{2, 1} = \frac 1m ((\Omega_\lambda)_{11}-(\Omega_\lambda)_{22})(\Omega^*_\lambda)_{12} - ((\Omega^*_\lambda)_{11} 
 - (\Omega^*_\lambda)_{22})(\Omega_\lambda)_{12} )
 \end{equation}
 and
 \begin{equation}\label{gauss-xi}
(\mathcal K -1) = \frac 1{m^2} \text{Tr}\Omega_\lambda\Omega^*_\lambda,
\end{equation}
where $\mathcal K$ is the Gaussian curvature of the M\"{o}bius metric $I^\xi = m|du|^2$. Of course, as one may verify, \eqref{codazzi-y}, \eqref{codazzi-y-star}, \eqref{codazzi-mix} and
\eqref{gauss-xi} are exactly the Gauss-Codazzi equations for the surface $\xi$ in the de Sitter spacetime $\mathbb{S}^{1, 3}$ induced by the conformal Gauss map
$\xi$ of the surface $\hat x$ in conformal 3-sphere $\mathbb{S}^3$.
\\

Now we are ready to state and prove the local fundamental theorem for surfaces in conformal round 3-sphere $\mathbb{S}^3$.

\begin{theorem}\label{local-fundamental} Suppose that, on a domain in $D\subset\mathbb{R}^2$,  we are given the following
\begin{itemize} 
\item a traceless symmetric 2-form $\Omega$
\item a positive function $m$ or equivalently $E$ such that $m= \frac {-\det\Omega}E$
\item a 1-form $\omega$
\item a symmetric 2-form $\Omega^*$.
\end{itemize}
And suppose that they satisfy the integrability conditions \eqref{codazzi-y} - \eqref{gauss-xi}. Then, for a given point $p_0$
in $D$, there exists an open neighborhood $D_0$ of $p_0$ in $D$,  a parametrized surface $\hat x: D_0\to\mathbb{S}^3$ 
with no umbilical point,  and a positive function $\hat\lambda: D_0\to \mathbb{R}^+$ with $\hat\lambda(p_0) = 1$, such that 
\begin{itemize}
\item $\Omega = \hat\lambda \overset{\circ}{II}$, where $\overset{\circ}{II}$ is the traceless part of the second fundamental form of $\hat x$ 
in the standard round $\mathbb{S}^3$
\item $m|du|^2= <d\xi, d\xi>$ is the M\"{o}bius metric induced by the conformal Gauss map $\xi$ of $\hat x$ 
\item $\omega  = -I((\overset{\circ}{II})^{-1}(d H)) - d \log\hat\lambda$, where $I$ is the first fundamental form and $H$ is the mean curvature of $\hat x$ 
in the standard round $\mathbb{S}^3$
\item $\Omega^* = -  \hat\lambda^{-1} <d\xi, dy^*>$, where $y^* = \frac 1{1-\hat x\cdot \hat x^*}(1, \hat x^*)$ and $\hat x^*$ is the conformal transform of $\hat x$.
\end{itemize}
The surface $\hat x$ is unique up to a conformal transformation of $\mathbb{S}^3$. 
\end{theorem}  

\proof We start with choosing starting values for $y, y^*, \xi_{u^1}, \xi_{u^2}, \xi$ at $p_0= (u_0^1, u_0^2)$, First we take a null vector 
$$
y(u^1_0, u^2_0) = y_0 = (1, \hat x_0)
$$
for some $\hat x_0\in \mathbb{S}^3\subset \mathbb{R}^4$. Then we choose $\xi(u^1_0, u^2_0)= \xi_0\in\mathbb{R}^{1,4}$ such that
\begin{equation}\label{initial-1}
<y_0, \xi_0> = 0 \text{ and } <\xi_0, \xi_0> = 1.
\end{equation}
Next we choose $\xi_{u^1}(u^1_0, u^2_0) = \xi^1_0\in\mathbb{R}^{1,4}$ and $\xi_{u^2}(u^1_0, u^2_0) = \xi^2_0\in\mathbb{R}^{1,4}$ 
such that
\begin{equation}\label{initial-2}
\aligned
<\xi^1_0, \xi^1_0>& =<\xi^2_0, \xi^2_0> = m(u^1_0, u^2_0), \\
<\xi^1_0, \xi^2_0> & = <\xi_0, \xi_0^1>=<\xi_0, \xi_0^2> = <y_0, \xi_0^1>=<y_0, \xi_0^2>=0.\endaligned
\end{equation}
Finally choose the unique null vector $y^*(u^1_0, u^2_0) = y^*_0$ such that
\begin{equation}\label{initial-3}
\aligned
<y^*_0, y_0> & = -1\\
<y^*_0, y^*_0> & = <y^*_0, \xi_0> = <y^*_0, \xi_0^1> = <y^*_0, \xi_0^2> = 0.\endaligned
\end{equation}
Notice that for any other choice of $\{y_1, y^*_1, \xi_1^1, \xi^2_1, \xi_1\}$ satisfying the same orthonormal properties in
\eqref{initial-1} - \eqref{initial-3}, there is a Lorentz transformation that takes one to the other. With the integrability conditions
assumed we may solve the systems \eqref{structure-equation-1} and \eqref{structure-equation-2} at least in an open neighborhood
$D_0$ of $p_0$ in $D$. Using the uniqueness of solutions to systems of linear ODE one sees that the solution $\{y, y^*, \frac 1{\sqrt m}\xi_{u^1}, \frac 1{\sqrt m}
\xi_{u^2}, \xi\}$ remains to be orthonormal in the Minkowski metric in $D_0$.\\

Now one should realize that the $y = \hat\lambda (1, \hat x)$ here is with some positive $\hat\lambda$ (not necessarily identically $1$ in $D_0$). It is then clear
from all previous calculations that the rest of the statements in the theorem can be easily verified.
\endproof

\subsection{The geometry of the associate ruled surface $x^+$ in hyperbolic space $\mathbb{H}^4$}

In this section we want to discuss the geometry of the associate ruled 3-surface $x^+$ in $\mathbb{H}^{4}$, which is associated with a given 
surface $\hat x$ in the conformal 3-sphere. It's relation to the associate 
surface $\tilde x$ is very much analogous to the one between the ambient spacetime and the Poincar\'{e}-Einstein manifold of a given conformal 
manifold in the work of Fefferman
and Graham. It deems to be useful to understand the geometry of the associate ruled 3-surface $x^+$ in $\mathbb{H}^4$.  \\

It is rather easy now to do calculations for $x^+$ after we have calculated the first fundamental form for the associate 4-surface $\tilde x$ in Minkowski spacetime 
$\mathbb{R}^{1,4}$ in section \ref{Sect:I-tilde}. We first have
$$
dx^+ = \frac 1{\sqrt 2} (e^ty_\lambda -e^{-t}y^*_\lambda)dt + (e^t(y_\lambda)_{u^1}+e^{-t}(y^*_\lambda)_{u^1})du^1 + (e^t(y_\lambda)_{u^2}+e^{-t}(y^*_\lambda)_{u^2})du^2
$$
and, using \eqref{tilde-x-I},
\begin{equation}\label{x-+-I}
\aligned
I^{x^+} & = (dt)^2 - 2\omega^\lambda_idtdu^i+(\frac {e^{2t}} {2m}((\Omega_\lambda)_{i1}(\Omega_\lambda)_{j1}+(\Omega_\lambda)_{i2}(\Omega_\lambda)_{j2}) \\
& \quad\quad +   (\omega_i\omega_j +  \frac 1m((\Omega_\lambda)_{i1}(\Omega^*_\lambda)_{j1}+(\Omega_\lambda)_{i2}(\Omega^*_\lambda)_{j2})) \\
& \quad\quad  + \frac {e^{-2t}}{2m}  ((\Omega^*_\lambda)_{i1}(\Omega^*_\lambda)_{j1}+(\Omega^*_\lambda)_{i2}(\Omega^*_\lambda)_{j2}))du^idu^j\endaligned
\end{equation}
One can calculate the determinant
\begin{equation}\label{+-determinant}
\det I^{x^+} = \frac 1{8m^2}(E^2_\lambda|e^t\Omega_\lambda + e^{-t}\Omega^*_\lambda|^2 - e^{-2t}((\Omega^*_\lambda)_{11} + (\Omega^*_\lambda)_{22})^2)^2,
\end{equation}
which can tell us where the associate ruled surface $x^+$ is degenerate.
\\

To obtain the second fundamental form of the surface $x^+$ it suffices to see that the conformal Gauss map $\xi$ is still the unit normal vector to the surface $x^+$ in the 
hyperboloid $\mathbb{H}^4$. Hence
\begin{equation}\label{x-+-II}
II^{x^+} = - <dx^+, d\xi> = \frac 1{\sqrt 2}(e^t\Omega_\lambda + e^{-t}\Omega^*_\lambda).
\end{equation}

By the similar calculations as that in the previous section we have  the mean curvature of the associate ruled surface $x^+$ as follows:
\begin{equation}\label{plus-H}
H^{x^+} = e^{-3t}\frac {\sqrt{2}\det\Omega_\lambda {\mathcal H}_\lambda}
{(\det\Omega_\lambda -e^{-2t}\text{Tr}\Omega_\lambda\Omega^*_\lambda + e^{-4t}\det\Omega^*_\lambda)}.
\end{equation}

\begin{theorem}\label{Th:+-H} Suppose that $\hat x$ is an immersed surface in the conformal sphere $\mathbb{S}^3$ with no umbilical point and that $x^+$ is the associate ruled surface in
the hyperboloid $\mathbb{H}^4$. Then $\hat x$ is a Willmore surface in the conformal sphere if and only if the associate ruled 3-surface $x^+$ in the hyperboloid is a minimal 
surface.
\end{theorem}

\section{Scalar invariants of surfaces in conformal round 3-sphere}

In this section we want to introduce scalar local invariants for surfaces in conformal round 3-sphere $\mathbb{S}^3$. We will first recall what are scalar 
invariants for hypersurfaces in (pseudo-)Riemannian geometry. Inspired by the work of Fefferman and Graham on scalar local invariants in conformal geometry 
we are going to use the associate surface $\tilde x$ in the Minkowski $\mathbb{R}^{1,4}$  of a given surface $\hat x$ in 3-sphere $\mathbb{S}^3$, where one considers
the standard conformal 3-sphere as the projectivized positive light cone of the Minkowski spacetime to construct scalar local invariant. 

\subsection{Scalar invariants of 4-surfaces in $\mathbb{R}^{1,4}$}

For our purpose we will focus on the discussion of scalar (pseudo-)Riemannian invariants of 4-surfaces $\tilde x$ in the Minkowski spacetime $\mathbb{R}^{1,4}$.
Suppose that 
$$
\phi = \phi(v^2, v^3, v^4, v^5): A\subset \mathbb{R}^4\to \mathbb{R}^{1,4}
$$
is a local parametrization of a surface $\tilde x$, where $A$ is a domain in $\mathbb{R}^4$. Hence it induces a local coordinate 
$$
\tilde \phi = \tilde \phi(v^1, v^2, v^3, v^4, v^5) : B\subset (-\epsilon, \epsilon)\times A\to \mathbb{R}^{1,4}
$$ 
for $\mathbb{R}^{1,4} $ such that
$$
\phi (v^2, v^3, v^4, v^5)  = \tilde \phi(0, v^2, v^3, v^4, v^5).
$$
We will use the Capital Latin letters to stand for indices from $1$ to $5$ and Latin letters to stand for the indices from $2$ to $5$.  And we will use $v = (v^1, v^2, \cdots, v^5)$
and $\hat v = (v^2, \cdots, v^5)$. Hence the Minkowski metric in this coordinate is give as
$$
\tilde{\mathcal{G}}_0 = <d\tilde\phi, d\tilde\phi> = (\tilde{\mathcal G}_0)_{IJ}dv^I dv^J
$$
and the fist fundamental form for $\tilde x$ in $\mathbb{R}^{1, 4}$ is given as
$$
I^{\tilde x} = <d\phi, d\phi> = \tilde g_{ij}dv^idv^j = (\tilde{\mathcal G}_0)_{ij}|_{v^1=0}dv^idv^j.
$$
To be more restrictive we will assume that the surface $\tilde x$ is timelike and let 
$$
\xi: B \to \mathbb{S}^{1,3}
$$
be a unit normal vector field on $\tilde x$ in $\mathbb{R}^{1,4}$. Then the second fundamental form for $\tilde x$ is given as

$$
II^{\tilde x} = - <d\phi, d\xi> = \tilde h_{ij}dv^idv^j.
$$
And
$$
\xi_{v^i} = - \tilde h_{ik}\tilde g^{kj} \phi_{v^j}.
$$
\begin{definition} Let ${\bf i}: \textup{M}^{n-1}\to \textup{N}^n$ be an immersed hypersurface and let $g$ be a (pseudo)-Riemannain metric on the ambient manifold $\textup{N}^n$.
A scalar (pseudo-)Riemannain invariant  $\textup{I}({\bf i}, \textup{N}^n, g)$ for the hypersurface ${\bf i}$ in $\textup{N}^n$ at a point $p_0$ on the surface ${\bf i}$
is a polynomial in the variables that are the coordinate partial derivatives of $g_{IJ}$ of any order and the reciprocal of the
determinant of $g_{IJ}$ at the point $p_0$ such that the value of  $\textup{I}({\bf i}, \textup{N}^n, g)$ at $p_0$ is independent of choices of local coordinates 
$\tilde\phi$ of $\textup{N}^n$ which are induced from a parametrization $\phi$ of the surface ${\bf i}$ nearby the given point $p_0$. 
\end{definition}

The well-known examples of scalar Riemannian invariants for $\tilde x$ in $\mathbb{R}^{1,4}$ are
\\

\begin{itemize}
\item $\tilde H = \tilde g^{ij}\tilde h_{ij}$
\item $|\tilde h|^2 = \tilde g^{ik}\tilde g^{jl}\tilde h_{ij}\tilde h_{kl}$ and $\tilde H^2 =\tilde g^{ij}\tilde g^{kl}\tilde h_{ij}\tilde h_{kl}$
\item $\tilde \Delta \tilde H = \tilde g^{kl} \tilde g^{ij}\tilde h_{ij,kl}$, $\text{Div}\text{Div}\tilde h = \tilde g^{ik}\tilde g^{jl}\tilde h_{ij, kl}$, 
$\tilde H|\tilde h|^2 = \tilde g^{ik}\tilde g^{jl}\tilde g^{mn}\tilde h_{ij}\tilde h_{kl}\tilde h_{mn}$, 
\newline 
$\quad \text{Tr}_{\tilde g}\tilde h^3 = \tilde g^{in}\tilde g^{jk}\tilde g^{km}\tilde h_{ij}\tilde h_{kl}\tilde h_{mn}$,  
and $\tilde H^3 =  \tilde g^{ij}\tilde g^{kl}\tilde g^{mn}\tilde h_{ij}\tilde h_{kl}\tilde h_{mn}$
\item $|\tilde\nabla \tilde h|^2 = \tilde g^{ip}\tilde g^{jq}\tilde g^{kr}\tilde h_{ij,k}\tilde h_{pq,r}$, $\ \tilde g^{ip}\tilde g^{jr}\tilde g^{kq}\tilde h_{ij,k}\tilde h_{pq,r}$, 
$\ \tilde g^{ip}\tilde g^{jr}\tilde g^{kq}\tilde h_{ij,k}\tilde h_{pq,r}$
\newline $|\tilde\nabla \tilde H|^2= \tilde g^{ij}\tilde g^{pq}\tilde g^{kr}\tilde h_{ij,k}\tilde h_{pq,r}$, 
$\ |\widetilde{\text{Div}}\tilde h|^2 = \tilde g^{ip}\tilde g^{jk}\tilde g^{qr}\tilde h_{ij,k}\tilde h_{pq,r}$, $\ \widetilde{\text{Div}}\tilde h \cdot d \tilde H$
\item $\tilde\Delta\tilde\Delta\tilde H$
\end{itemize}

\vskip 0.1in
Each scalar invariant has an order. To find the order of each scalar invariant one simply scales the metric by a constant $\kappa$ and see 
what is the dimension of the scalar invariant. For example, we can easily find that

$$
\aligned
\tilde H [\kappa^2 \tilde {\mathcal G}_0] & = \kappa^{-1}\tilde H [\tilde{\mathcal G}_0]\\
|\tilde h|^2 [\kappa^2 \tilde {\mathcal G}_0] & = \kappa^{-2} |\tilde h|^2[\tilde {\mathcal G}_0]\\
\tilde \Delta \tilde H[\kappa^2 \tilde {\mathcal G}_0] & = \kappa^{-3} \tilde\Delta\tilde H [\tilde {\mathcal G}_0]\\
|\tilde\nabla\tilde h|^2 [\kappa^2 \tilde {\mathcal G}_0] & =\kappa^{-4} |\tilde\nabla\tilde h|^2[\tilde {\mathcal G}_0]\\
\tilde\Delta\tilde\Delta\tilde H [\kappa^2\tilde {\mathcal G}_0] & = \kappa^{-5}\tilde\Delta\tilde\Delta\tilde H [\tilde {\mathcal G}_0].
\endaligned
$$
\\

To understand what are scalar Riemannian invariants $\textup{I}(\tilde x, \mathbb{R}^{1,4}, \tilde{\mathcal G}_0)$ 
we want to use the so-called Fermi coordinates. A Fermi coordinate is one such that
1) on the surface $\phi$ is a normal coordinate at a given point $\tilde x_0$; 2) the coordinate curves $\tilde\phi(t, v^2, v^3, v^4, v^5)$ 
is a geodesic perpendicular to the surface at $\phi(v^2, v^3, v^4, v^5)$ with unit speed (a line segment perpendicular to the surface in $\mathbb{R}^{1, 4}$). Hence, for a
Fermi coordinate,
\begin{equation}\label{line-segment}
\tilde\phi (v^1, \cdots, v^5) = \phi(v^2, \cdots, v^5) + v^1\xi.
\end{equation}
The following facts are well known.

\begin{lemma} Suppose that $\tilde x$ is a timelike hypersurface in $\mathbb{R}^{1,4}$. Suppose that $\tilde\phi$ is a Fermi coordinate at a given point $\tilde x_0$. 
Then
$$
\tilde{\mathcal G}_0 = \left[\begin{matrix} 1 & 0 \\0 & [{\mathcal G}_{ij}]\end{matrix}\right] 
$$
and 
$$
{\mathcal G}_{ij}(v^1, \hat v) = \tilde g_{ij}(\hat v) - 2\tilde h_{ij}(\hat v) v^1 + \tilde h_{ik}(\hat v)\tilde h_{jl}(\hat v)\tilde g^{kl}(\hat v) (v^1)^2,
$$
where
$$
\tilde g_{ij} (\hat v) = \eta_{ij} -\frac 23 \tilde R_{ikjl}v^kv^l + \cdots
$$
$$
\tilde h_{ij} (\hat v) = \tilde h_{ij}(0) + \tilde h_{ij, k}(0)v^k + \cdots
$$
$\tilde R_{ijkl} = \tilde h_{ik}\tilde h_{jl} - \tilde h_{ij}\tilde h_{kl}$ is the Riemann curvature tensor for $\tilde x$
and $\eta$ is standard matrix of signature $\{-1, 1, 1, 1\}$. Moreover all the coefficients
in the Taylor's expansions for $G_{ij}$  are polynomials of $\tilde h_{ij}$ and the covariant 
derivatives of $\tilde h_{ij}$ at $\tilde x_0$.
\end{lemma}

Therefore, in the light of Weyl theorem on the invariants of orthogonal groups,  we may conclude that

\begin{proposition} All scalar invariants $\textup{I}(\tilde x, \mathbb{R}^{1,4}, \tilde{\mathcal G}_0)$ 
of a surface $\tilde x$ in $\mathbb{R}^{1,4}$ are linear combinations of terms that are complete contractions of tensor product of 
the second fundamental form $\tilde h$ and the covariant derivatives of $\tilde h$.
\end{proposition}
\proof From the above lemma it is easily that all scalar invariants of a surface $\tilde x$ in $\mathbb{R}^{1, 4}$ are polynomials of the first fundamental form $\tilde g$, 
the second fundamental form $\tilde h$ and covariant derivatives of the second fundamental form $\tilde h$, if we evaluate them in a Fermi coordinate for the surface.
Then, by the Weyl theorem on the invariants of orthogonal groups, we know they are linear combinations of full contractions of $\tilde h$ and covariant derivatives of $\tilde h$.
\endproof

\subsection{Scalar invariants of the homogeneous associate surface $\tilde x$ in $\mathbb{R}^{1,4}$}\label{calculations}

Let us work with the parametrization 
$$
\tilde x = \alpha \hat\lambda(1, \hat x) + \alpha\rho \hat\lambda^{-1} \frac 1{1-a} (1, \hat x^*) = \alpha y_\lambda + \alpha\rho y^*_\lambda
$$
and use the calculations given in Section \ref{Sect:I-tilde} and Section \ref{Sect:II-tilde}.
Now let us compute some scalar invariants for our associate surface $\tilde x$ on the light cone where $\rho = 0$. Then the first fundamental form is
$$
I_{\tilde x} |_{\rho = 0}= \left[\begin{matrix} \begin{matrix} 0 & -\alpha \\-\alpha & 0 \end{matrix} &   
\begin{matrix} 0 & 0 \\ \alpha^2\omega_1^\lambda & \alpha^2 \omega_2^\lambda \end{matrix} \\
 \begin{matrix} 0 &\quad\alpha^2\omega_1^\lambda \\ 0  & \quad \alpha^2\omega_2^\lambda \end{matrix} & 
  \begin{matrix} \alpha^2 E_\lambda & 0 \\ 0  & \alpha^2 E_\lambda \end{matrix} \end{matrix}\right]
 $$
from \eqref{matrix-tilde},  whose inverse is
 $$
 I_{\tilde x}^{-1}|_{\rho = 0} = \left[\begin{matrix} \begin{matrix} |\omega^\lambda|^2 & -\frac 1\alpha \\-\frac 1\alpha & 0 \end{matrix} &   \begin{matrix} 
 \frac {\omega^\lambda_1}{\alpha E_\lambda}
  & \frac {\omega^\lambda_2}{\alpha E_\lambda} \\ 0 & 0 \end{matrix} \\
 \begin{matrix} \frac {\omega^\lambda_1}{\alpha E_\lambda} & \quad 0 \\ 
 \frac {\omega^\lambda_2}{\alpha E_\lambda}  & \quad 0 \end{matrix} &  \begin{matrix} \frac 1{\alpha^2 E_\lambda} & 0 \\ 0  & 
 \frac 1{\alpha^2 E_\lambda} \end{matrix} \end{matrix}\right].
 $$
And the second fundamental form at $\rho=0$ is
 $$
 II_{\tilde x}|_{\rho = 0} = \left[\begin{matrix} 0 & 0 \\
 0 & \alpha\Omega_\lambda\end{matrix}\right].
 $$
% Therefore
% \begin{equation}\label{2-form}
 %I_{\tilde x}^{-1}II_{\tilde x}|_{\rho = 0} = \left[\begin{matrix} \begin{matrix} 0  & 0 \\ 0 & 0 \end{matrix} &   \begin{matrix} -H_{u^1} & -H_{u^2} \\0 & 0 \end{matrix} \\
% \begin{matrix} 0 & 0 \\ 0  & 0 \end{matrix} &  \begin{matrix} \frac 1{\alpha E}\Omega_{11} & \frac 1{\alpha E}\Omega_{12} \\ \frac 1{\alpha E}\Omega_{21} 
 % & \frac 1{\alpha E}\Omega_{22} \end{matrix} \end{matrix}\right]
 %\end{equation}

So the simplest (pseudo-)Riemannian invariants is the mean curvature $\tilde H$, but it is clear that 
$$
\tilde H |_{\rho = 0} = \frac 1{\alpha E_\lambda} ((\Omega_\lambda)_{11} + (\Omega_\lambda)_{22}) = 0.
$$ 
The first non-trivial one is 
\begin{equation}\label{1st-one}
|\tilde h|^2|_{\rho =0} = \tilde g^{ik}\tilde g^{jl}\tilde h_{ij}\tilde h_{kl} |_{\rho = 0} = \alpha^{-2}|\Omega_\lambda|^2,
\end{equation}
which produces the first non-trivial invariant $|\overset{\circ}{II}|^2$ for the surface $\hat x$ in the conformal 3-sphere(cf. see the definition for scalar 
invariant of surfaces in the conformal 3-sphere in the next subsection). In fact the following non-trivial invariants without taking any derivative are
all easy to calculate
$$
\text{Tr}_{I^{\tilde x}}\tilde h^k|_{\rho = 0}  = \alpha^{-k} \text{Tr}_{I^{\hat x}_\lambda} \Omega_\lambda^k
$$
for any $k = 2, 3, \cdots.$ Obviously those are the ones that can been easily seen with no difficulty at all. 
\\

Next we want to calculate $|\nabla \tilde H|^2$ and $\tilde\Delta\tilde H$ at $\rho = 0$. To do so, let us first recall from Section \ref{Sect:II-tilde} the mean curvature
$$
\tilde H = \frac {\rho \det\Omega_\lambda{\mathcal H}_\lambda}
{\alpha(\det\Omega_\lambda - \rho\text{Tr}\Omega_\lambda\Omega^*_\lambda + \rho^2\det\Omega^*_\lambda)}.
$$
Hence $\tilde H_\alpha = \tilde H_{u^1} = \tilde H_{u^2} = 0$ and $|\nabla \tilde H|^2 = 0$ at $\rho=0$, that is, $|\nabla \tilde H|^2$ gives no invariant for the surface $\hat x$. 
Let us set the convention to have $a, b, c$ stand for $\alpha, \rho$; $i,j,k$ stand for $u^1, u^2$, and $A, B, C$ stand for all four variables. 
We then calculate, at $\rho =0$, 
\begin{equation}\label{laplace-h}
\aligned
\tilde\Delta\tilde H & = \frac 1{\sqrt{|\tilde g|}} \partial_A (\sqrt {|\tilde g|} \tilde g^{AB}\partial_B\tilde H)\\
& =\frac 1{\alpha^3 E} (\partial_\alpha (\sqrt {|\tilde g|} \tilde g^{\alpha \rho}
\partial_\rho\tilde H) + \partial_\rho (\sqrt {|\tilde g|} \tilde g^{\rho B}\partial_B\tilde H) + \partial_i(\sqrt {|\tilde g|} \tilde g^{i\rho}\partial_\rho\tilde H))\\
& =\frac 1{\alpha^3E}(\partial_\alpha (\sqrt {|\tilde g|} \tilde g^{\alpha \rho}
\partial_\rho\tilde H) +  \partial_\rho (\sqrt{|\tilde g|} g^{\rho\alpha}\partial_\alpha\tilde H) + \sqrt{|\tilde g|} (\partial_\rho\tilde g^{\rho\rho}) \partial_\rho \tilde H)\\
& = 2\alpha^{-3}{\mathcal H}_\lambda
\endaligned
\end{equation}
where one needs to use the fact that $\tilde g^{\rho\rho}|_{\rho = 0} = 0$ and $\partial_\rho \tilde g^{\rho\rho}|_{\rho = 0} = \frac {2}{\alpha^2}$ based on calculations \eqref{app-2j}
in Appendix \ref{app-inverse-tilde}. This confirms that ${\mathcal H}_\lambda$ is indeed a conformal invariant of order $3$ for a surface $\hat x$ in 3-sphere in 
general conformal metric $\lambda^2g_0$.
\\

The next invariant we want to calculate is $\tilde\Delta\tilde\Delta\tilde H$. To do so, from \eqref{app-2j} in Appendix \ref{app-inverse-tilde}, we observe the following:
\begin{equation}\label{der-inverse}
\aligned
\partial_\rho|_{\rho=0}\tilde g^{\rho\alpha} & = -\frac 2\alpha |\omega^\lambda|^2, \quad \partial_\rho|_{\rho=0}\tilde g^{\rho\rho} = \frac 2{\alpha^2}, \\
\partial_\rho|_{\rho=0}\tilde g^{\rho i} & = -\frac 2{\alpha^2} \frac {\omega^\lambda_i}{E_\lambda} \text{ and } \partial_\rho\partial_\rho|_{\rho=0}\tilde g^{\rho\rho}
= \frac 8{\alpha^2}|\omega^\lambda|^2.\endaligned
\end{equation}
After a lengthy calculation we get
\begin{equation}\label{double-laplace}
\aligned
\tilde\Delta & \tilde\Delta\tilde H|_{\rho=0} = 8\alpha^{-5} (\Delta_\lambda{\mathcal H}_\lambda +9 |\omega^\lambda|^2{\mathcal H}_\lambda - 3 \text{Div}({\omega^\lambda})
{\mathcal H}_\lambda \\ & \quad - 6\omega^\lambda(\nabla{\mathcal H}_\lambda) 
- \frac {3\text{Tr}(\Omega_\lambda\Omega^*_\lambda)}{2m^2} |\Omega_\lambda|^2{\mathcal H}_\lambda).\endaligned
\end{equation}
This tells us that $\Delta_\lambda{\mathcal H}_\lambda +9 |\omega^\lambda|^2{\mathcal H}_\lambda - 3 \text{Div}({\omega^\lambda})
{\mathcal H}_\lambda - 6\omega^\lambda(\nabla{\mathcal H}_\lambda)- \frac {3\text{Tr}(\Omega_\lambda\Omega^*_\lambda)}{2m^2} |\Omega_\lambda|^2{\mathcal H}_\lambda$ is a conformal invariant of order $5$ for the surface $\hat x$ in 3-sphere.
\\

We can also calculate the covariant derivatives of the second fundamental forms for the associate surface. We first list the relevant Christoffel symbols for the calculation
\begin{equation}\label{christoffel}
\aligned
\tilde\Gamma^k_{\alpha\alpha} & = \tilde\Gamma^k_{\rho\rho} = \tilde\Gamma^k_{\alpha\rho} = 0 \\
\tilde\Gamma^k_{\alpha j} & = \alpha^{-1} \delta_{jk} \\ 
\tilde\Gamma^k_{\rho j}  & = \frac 1{2E_\lambda} ((\omega^\lambda_k)_{u^j} - (\omega^\lambda_j)_{u^k} + \frac 1m((\Omega_\lambda)_{jl}(\Omega^*_\lambda)_{kl}
+(\Omega_\lambda)_{kl}(\Omega^*_\lambda)_{jl})) \\
%- 2(\Omega_\lambda)_{kj} \frac {{\mathcal H}_\lambda}{|\Omega_\lambda|^2}) + \frac {\text{Tr}(\Omega_\lambda\Omega^*_\lambda)}{4m^2}|\Omega_\lambda|^2\delta_{jk} \\
\tilde\Gamma^k_{ij} & = (\Gamma_\lambda)^k_{ij} - \omega^\lambda_k\delta_{ij}.
\endaligned
\end{equation}
Then we calculate
\begin{equation}\label{co-derivative}
\aligned
\tilde h_{ab, C} &= 0\\
\tilde h_{ai,b} & = 0 \\
\tilde h_{\alpha j, k} & = -(\Omega_\lambda)_{jk} \\
\tilde h_{\rho j,k}  & = - \frac \alpha{2E_\lambda}((\Omega_\lambda)_{ij}((\omega^\lambda_i)_{u^k} - (\omega^\lambda_k)_{u^i} +\frac 1m ((\Omega_\lambda)_{kl}(\Omega_\lambda^*)_{il}
+ (\Omega_\lambda)_{il}(\Omega^*_\lambda)_{kl}) \\
\tilde h_{ij,\alpha} & = -(\Omega_\lambda)_{ij} \\
\tilde h_{ij, \rho}  & =\alpha(\Omega_\lambda^*)_{ij} \\ & \quad\quad -  \frac \alpha{2 E}((\Omega_\lambda)_{lj} ((\omega^\lambda_l)_{u^i}  - (\omega^\lambda_i)_{u^l} +\frac 1m 
((\Omega_\lambda)_{kl}(\Omega_\lambda^*)_{ki} + (\Omega_\lambda)_{ki}(\Omega^*_\lambda)_{kl})\\
& \quad\quad -\frac \alpha{2E}(\Omega_\lambda)_{il}((\omega^\lambda_l)_{u^j} - (\omega^\lambda_j)_{u^l} +\frac 1m 
((\Omega_\lambda)_{kl}(\Omega_\lambda^*)_{kj} + (\Omega_\lambda)_{kj}(\Omega^*_\lambda)_{kl})\\
\tilde h_{ij, k} & = \alpha(\Omega_\lambda)_{ij,k} + \alpha(\Omega_\lambda)_{lj}\omega^\lambda_l\delta_{ik} + \alpha(\Omega_\lambda)_{il}\omega^\lambda_l\delta_{jk}.
\endaligned 
\end{equation}
The easy one is 
 $$
\phi_\alpha =\tilde h_{\alpha j,k}\tilde g^{jk} = 0 \text{ and } \phi_\rho = \tilde h_{\rho j,k}\tilde g^{jk} = \frac 1\alpha{\mathcal H}_\lambda
 $$
 in the light of \eqref{trace-star}. While
 $$
 \aligned
 \phi_i  & = \tilde h_{iB,C}\tilde g^{BC}  = \tilde h_{ij,C}\tilde g^{jC} + \tilde h_{ib, k}\tilde g^{bk}
 = \tilde h_{ij,k}\tilde g^{jk} + \tilde h_{ij,\alpha}\tilde g^{j\alpha}  + \tilde h_{i\alpha,k}\tilde g^{\alpha k}\\
 & = \frac 1{\alpha E_\lambda}(\Omega_\lambda)_{ij,j}  +  \frac 3{\alpha E_\lambda}(\Omega_\lambda)_{ij}\omega^\lambda_j 
 - \frac 1{\alpha E_\lambda}(\Omega_\lambda)_{ij}\omega^\lambda_j  -\frac 1{\alpha E_\lambda}(\Omega_\lambda)_{ij}\omega^\lambda_j\\
 & = \frac 1{\alpha E_\lambda}(\Omega_\lambda)_{ij,j}  +  \frac 1{\alpha E_\lambda}(\Omega_\lambda)_{ij}\omega^\lambda_j = 0
 \endaligned
 $$
 due to the integrability condition \eqref{codazzi-y}. Thus $|\widetilde{\text{Div}}\tilde h|^2(=0)$ does not give any invariant on the surface $\hat x$, nor does
 $\widetilde{\text{Div}}\tilde h\cdot d\tilde H (= 0)$. Because $\tilde g^{\rho\rho} |_{\rho=0}= 0$. 
 \\

We want to calculate $|\tilde\nabla \tilde h|^2$ since we have all the covariant derivatives $\tilde h_{AB,C}$ in \eqref{co-derivative}. The calculation is direct yet very long.
We omit the detail here.
\begin{equation*}
|\tilde\nabla\tilde h|^2 |_{\rho = 0} =  \alpha^{-4}(|\nabla\Omega|^2 +8 |d H|^2 -6\Omega\cdot\Omega^* - \frac 2{E_\lambda^3}(\Omega_\lambda)_{ij}\omega^\lambda_k(R^\lambda)_{3ijk}-
\frac 6{E_\lambda^3}(\Omega_\lambda)_{ij}(\Omega_\lambda)_{ki,j}\omega^\lambda_k),
\end{equation*}
where the Codazzi equation for the surface $\hat x$ in $(\mathbb{S}^3, \ \lambda^2g_0)$
$$
(\Omega_\lambda)_{ij,k} = (\Omega_\lambda)_{ik,j}  + (R^\lambda)_{3ijk} + (H_\lambda)_{u^j}E_\lambda\delta_{ik} - (H_\lambda)_{u^k} E_\lambda\delta_{ij}
$$
has been used. At this point we like to write each term as local scalar invariant of the surface $\hat x$ in $(\mathbb{S}^3, \ \lambda^2g_0)$. 
We first calculate
$$
\aligned
& (\Omega_\lambda)_{ij} \omega^\lambda_k(R^\lambda)_{3ijk} = (\Omega_\lambda)_{ij}\omega^\lambda_1(R^\lambda)_{3ij1} + (\Omega_\lambda)_{ij}\omega^\lambda_2(R^\lambda)_{3ij2}\\
& =E_\lambda( (\Omega_\lambda)_{11}\omega^\lambda_1(R^\lambda)_{31} + (\Omega_\lambda)_{21}\omega^\lambda_1(R^\lambda)_{32} 
 + (\Omega_\lambda)_{22}\omega^\lambda_2(R^\lambda)_{32} + (\Omega_\lambda)_{12}\omega^\lambda_2(R^\lambda)_{31})\\
 & =  E_\lambda(\Omega_\lambda)_{ij} \omega^\lambda_j(R^\lambda)_{3i}  = -  E_\lambda^2 (H_\lambda)_{u^i}(R^\lambda)_{3i} = - 
 Ric^\lambda(\overset{\rightarrow}{\bf n}_\lambda, \lambda H_\lambda).
\endaligned
$$
Then we deal with the last term
$$
\aligned
(\Omega_\lambda)_{ij}  &  (\Omega_\lambda)_{ki,j}\omega^\lambda_k = (\Omega_\lambda)_{ij}((\Omega_\lambda)_{ki}\omega^\lambda_k)_{,j} - 
(\Omega_\lambda)_{ij}(\Omega_\lambda)_{ki}\omega^\lambda_{k,j}
\\ & = - E_\lambda(\Omega_\lambda)_{ij} (H_\lambda)_{i,j} - \frac 12 E_\lambda^3|\Omega_\lambda|^2\text{Div}(\omega^\lambda)
\endaligned
$$
where
$$
\aligned
\text{Div}(\omega^\lambda) & = E^{-1}_\lambda\omega^\lambda_{i,i} = E^{-1}_\lambda(\omega^\lambda_i)_{u^i} \\ & 
= E^{-1}_\lambda(<\Delta_0 y_\lambda, y^*_\lambda> + <(y_\lambda)_{u^i}, (y^*_\lambda)_{u^i}>)\\
& = H_\lambda^2 -  |\omega^\lambda|^2 + (R^\lambda)_{1212} +E^{-1}<(y_\lambda)_{u^i}, (y^*_\lambda)_{u^i}>\\
& = H^2_\lambda + 2 \frac {\Omega_\lambda\cdot\Omega^*_\lambda}{|\Omega_\lambda|^2} + E_\lambda^{-1}(R^\lambda)_{1212}\\
\Delta_0 y_\lambda & = 2E_\lambda H_\lambda \overset{\rightarrow}{\bf n}_\lambda + 2 E_\lambda y^\dagger_\lambda  - (R^\lambda)_{1212}y_\lambda
\endaligned
$$
and
$$
\sum_{i=1}^2 E^{-1}<(y)_{u^i}, (y^*)_{u^i}> = |\omega^\lambda|^2 + 2 \frac {\Omega_\lambda\cdot\Omega^*_\lambda}{|\Omega_\lambda|^2}.
$$
So we have obtained
\begin{equation}\label{norm-co-der}
\aligned 
|\nabla\tilde h|^2|_{\rho =0} & =\alpha^{-4}( |\nabla\Omega_\lambda|^2 + 8|dH_\lambda|^2 + 2 Ric^\lambda(\overset{\rightarrow}{\bf n}_\lambda, \nabla H_\lambda) 
+ 3 H^2_\lambda|\Omega|^2 \\ & +3K^T_\lambda |\Omega_\lambda|^2+ 6 \Omega_\lambda\cdot \text{Hess} (H_\lambda)) \endaligned
\end{equation}
where
$$
K^T_\lambda = E_\lambda^{-1}(R^\lambda)_{1212}
$$ 
is the sectional curvature of $(\mathbb{S}^3, \lambda^2 g_0)$ of the tangent plane to the surface $\hat x$.
\\

\subsection{Scalar invariants for surfaces in the conformal round 3-sphere}

Let us start with the definition of scalar invariants for surfaces in conformal sphere. 

\begin{definition} Let ${\bf i}: \textup{M}^{n-1}: \textup{N}^n$ be an immersed hypersurface and let $[g]$ be a class of conformal metrics on the ambient manifold $\textup{N}|^n$.
 $\textup{I}_c({\bf i}, \textup{N}^n, g)$ is said to be a scalar conformal invariant of the hypersurface ${\bf i}$ in the conformal manifold $(\textup{N}^n, [g])$ 
 if it is a scalar Riemannian invariant and
\begin{equation}\label{conformal-invariance}
\textup{I}_c({\bf i}, \textup{N}^n, \lambda^2g) = \lambda^{-k} \textup{I}_c({\bf i}, \textup{N}^n, g).
\end{equation}  
for any positive function $\lambda$ on $\textup{N}^n$, where $k$ is the order of the invariant $\textup{I}_c({\bf i}, \textup{N}^n, g)$. 
\end{definition}

Recall that, for an immersed surface $\hat x$ in $\mathbb{S}^3$, we have
$$
\overset{\circ}{II} (\hat x, \mathbb{S}^3, \lambda^2 g_0) = \lambda \overset{\circ}{II} (\hat x, \mathbb{S}^3, g_0).
$$
Hence it is easy to observe that
$$
|\overset{\circ}{II}|^2(\hat x, \mathbb{S}^3, \lambda^2 g_0) = \lambda^{-2}g_0^{ik}\lambda^{-2}g_0^{jl} \lambda\overset{\circ}{II}_{ij}\lambda\overset{\circ}{II}_{kl} 
= \lambda^{-2}\|\overset{\circ}{II}\|^2 (\hat x, \mathbb{S}^3, g_0)
$$
and
$$
\text{Tr}_{\lambda^2 g_0}(\overset{\circ}{II})^k(\hat x, \mathbb{S}^3, \lambda^2 g_0) = \lambda^{-k} \text{Tr}_{g_0}(\overset{\circ}{II})^k(\hat x, \mathbb{S}^3, g_0) \text{ for all $k=2, 3, \cdots$}.
$$
On the other hand, it does not seem easy to directly verify that ${\mathcal  H}_\lambda$ is a 
conformal invariant for a surface in the conformal 3-sphere, though this is a well-known one.
We have verified this in computing the mean curvature (cf. \eqref{mean-curvature-xi}) of the surface $\xi$ in the de Sitter 
spacetime $\mathbb{S}^{1,3}$ as well as in the above calculation of $\tilde\Delta \tilde H$ (cf. \eqref{laplace-h}) of the homogeneous associate surface $\tilde x$. In general it takes 
tremendous, if not impossible, to verify whether an invariant $\textup{I}(\hat x, \mathbb{S}^3, \lambda^2g_0)$ is conformally invariant, complicated by the six integrability conditions. 
The most important application of the construction of associate homogeneous surfaces is the following:

\begin{theorem}\label{main-inv} Suppose that $\hat x:\textup{M}^2\to\mathbb{S}^3$ is an immersed surface with no umbilical point. And suppose that
$$
\tilde x = \alpha y  + \alpha\rho y^*: \mathbb{R}^+\times\mathbb{R}^+\times\textup{M}^2\to\mathbb{R}^{1,4}
$$
is the associate surface for $\hat x$, where $\hat x^*$ is the conformal transform of $\hat x$. Then any scalar (pseudo)-Riemannian 
invariant $\textup{I}(\tilde x, \mathbb{R}^{1,4}, \tilde{\mathcal G}_0)$ evaluated  at $\rho = 0$, if it is nontrivial, is a scalar conformal invariant 
$\textup{I}_c(\hat x, \mathbb{S}^3, \lambda^2g_0)$ multiplied with $|\overset{\circ}{II}_\lambda|^{2n}$ for some integer $n$.
\end{theorem}
\proof For any invariant $\textup{I}(\tilde x, \mathbb{R}^{1,4}, \tilde{\mathcal G}_0)$, we know that it is a full contraction of tensor product of 
the second fundamental form and the covariant derivatives. For a choice of representative $\lambda^2 g_0$ on $\mathbb{S}^3$, in the corresponding
parametrization \eqref{lambda-para}, we claim that
\begin{equation}\label{Phi-lambda}
\textup{I}(\tilde x, \mathbb{R}^{1,4}, \tilde{\mathcal G}_0)|_{\rho = 0} = \alpha^{-k} \textup{I}(\hat x, \mathbb{S}^3, 
\lambda^2 g_0)|\overset{\circ}{II}_\lambda|^{2n}
\end{equation}
for a positive integer $k$ and a nonnegative integer $n$, due to the homogeneity of the associate surface. To see
the right side of \eqref{Phi-lambda} is indeed a scalar Riemannian invariant multiplied with factor $|\overset{\circ}{II}_\lambda|^{2n}$ for some integer $n$, we consider 
the tensors that determines the first and second fundamental forms of the associate surface in that parametrization. We recall from \eqref{Omega} that
$$
\Omega_\lambda = \overset{\circ}{II}_\lambda
$$
is the traceless part of the second fundamental form for the surface $\hat x$ in the 3-sphere with the conformal metric $\lambda^2 g_0$. We also know from \eqref{omega} 
that
$$
\omega^\lambda = - I_\lambda(( \overset{\circ}{II}_\lambda)^{-1}(dH_\lambda)) =  - \frac  2{| \overset{\circ}{II}_\lambda|^2} \overset{\circ}{II}_\lambda (\nabla H_\lambda),
$$
which causes  us to include the possible negative $n$ in the right side of \eqref{Phi-lambda}. We may also recall from \eqref{m-lambda} that
$$
m = \frac 12 E_\lambda |\overset{\circ}{II}_\lambda|^2.
$$
Next we want to show that $\Omega^*_\lambda$ is also  a tensor product of covariant 
derivatives of the 1-form $\omega^\lambda$, covariant derivatives of the second fundamental form $II_\lambda$ and 
covariant derivatives of Riemann curvature tensor of the conformal metric $\lambda^2 g_0$ on the 3-sphere(including 0th order). Recall the definition 
$$
(\Omega^*_\lambda)_{ij}  = < y^*_\lambda, \xi_{u^iu^j}>.
$$
We use the same idea in the calculation of  the trace of $\Omega^*$ in Section \ref{geometry-xi}. Hence we write
\begin{equation}\label{xi-ij-proj}
\xi_{u^iu^j}  = -(\Omega^*_\lambda)_{ij}y_\lambda - (\Omega_\lambda)_{ij}y^*_\lambda + (\Gamma_m)^k_{ij}\xi_{u^k} - m\delta_{ij}\xi.
\end{equation}
From \eqref{y-star-lambda} we know that 
$$
<y^*_\lambda, y^\dagger_\lambda> =  -  \frac 12(|\omega^\lambda|^2 + H_\lambda^2).
$$
Using $\xi = H_\lambda y_\lambda + \overset{\rightarrow}{\bf n}_\lambda$ from Lemma \ref{Lem:xi-lambda} and \eqref{n-dagger}, we have
$$
<\xi_{u^k}, y^\dagger_\lambda> = - (H_\lambda)_{u^k} +(R^\lambda)_{3k}
$$
and 
$$
<\xi, y^\dagger_\lambda> = - H_\lambda.
$$ 
Therefore we derive from \eqref{xi-ij-proj} that 
\begin{equation}\label{Omega-star-first}
<\xi_{u^iu^j}, y^\dagger_\lambda> = (\Omega^*_\lambda)_{ij} + \frac 12(|\omega^\lambda|^2 + H^2_\lambda)(\Omega_\lambda)_{ij} +(\Gamma_m)^k_{ij}(-H_{u^k}
+ (R^\lambda)_{3k}) + Hm\delta_{ij},
\end{equation}
where 
$$
(\Gamma_m)^k_{ij} = \Gamma^k_{ij} + \frac 12 |\Omega_\lambda|^{-2}(|\Omega_\lambda|^2_{u^i}\delta_{jk} + |\Omega_\lambda|^2_{u^j}\delta_{ik} - 
|\Omega_\lambda|^2_{u^k}\delta_{ij})
$$
is the Christofel symbols for the M\"{o}bius metric $m|du|^2$. 
On the other hand we have
$$
\xi_{u^iu^j} = (H_\lambda)_{u^iu^j} y_\lambda + (H_\lambda)_{u^i}(y_\lambda)_{u^j} + (H_\lambda)_{u^j}(y_\lambda)_{u^i} + H_\lambda (y_\lambda)_{u^iu^j} +
(\overset{\rightarrow}{\bf n}_\lambda)_{u^iu^j}
$$
which implies
\begin{equation}\label{Omega-star-second}
\aligned
<\xi_{u^iu^j}, y^\dagger_\lambda> & = - (H_\lambda)_{u^iu^j} + H_\lambda<(y_\lambda)_{u^iu^j}, y^\dagger_\lambda> + < (\overset{\rightarrow}{\bf n}_\lambda)_{u^iu^j},
y^\dagger_\lambda>\\
& =- (H_\lambda)_{u^iu^j} - H_\lambda<(y_\lambda)_{u^i}, (y^\dagger_\lambda)_{u^j}> - < (\overset{\rightarrow}{\bf n}_\lambda)_{u^i}, (y^\dagger_\lambda)_{u^j}>\\
& \quad\quad\quad  - <\overset{\rightarrow}{\bf n}_\lambda, (y^\dagger_\lambda)_{u^i}>_{u^j}\\
& =- (H_\lambda)_{u^iu^j}+ \frac 1{E_\lambda} (\Omega_\lambda)_{ik}<(y_\lambda)_{u^k}, (y^\dagger_\lambda)_{u^j}> 
 - <\overset{\rightarrow}{\bf n}_\lambda, (y^\dagger_\lambda)_{u^i}>_{u^j}\\
& =- (H_\lambda)_{u^iu^j} - \frac 1{E_\lambda} (\Omega_\lambda)_{ik}(R^\lambda)_{i3k3}  +  ((R^\lambda)_{3i})_{u^j}.
\endaligned
\end{equation}
by \eqref{i-dagger-j} and \eqref{n-dagger}. Thus, comparing \eqref{Omega-star-first} and \eqref{Omega-star-second}, we have
\begin{equation}\label{Omega-star-ij}
\aligned
(\Omega^*_\lambda)_{ij} & =  - (H_\lambda)_{u^i,u^j}  - H_\lambda m\delta_{ij} - \frac 1{E_\lambda} (\Omega_\lambda)_{ik}(R^\lambda)_{j3k3}  +  ((R^\lambda)_{3i})_{,u^j} 
\\ & \quad - \frac 12(|\omega^\lambda|^2 + H^2_\lambda)(\Omega_\lambda)_{ij}   \\ & \quad
+ \frac 12 |\Omega_\lambda|^{-2}(|\Omega_\lambda|^2_{u^i}\delta_{jk} + |\Omega_\lambda|^2_{u^j}\delta_{ik} - |\Omega_\lambda|^2_{u^k}\delta_{ij})
((H_\lambda)_{u^k} - (R^\lambda)_{3k}).
\endaligned
\end{equation}
The last factor that goes into the left side of the equation \eqref{Phi-lambda} is the reciprocal of the determinant: 
$$
\det\tilde g|_{\rho =0} =  - \frac {\alpha^6}{m^2}(pr - q^2)^2 |_{\rho=0} = \frac {\alpha^6}{m^2} (\det\Omega_\lambda)^2 = \alpha^6 E^2_\lambda = \alpha^6 \det I^{\hat x}_\lambda.
$$
due to \eqref{det-I-tilde}, where $I^{\hat x}_\lambda = (\hat x)^*(\lambda^2 g_0) = E_\lambda |du|^2$. 
\\

To verify that the right side of \eqref{Phi-lambda} is actually a conformal invariant, for a positive functions $\lambda$ on 3-sphere, we simply compare the right side of \eqref{Phi-lambda} evaluated at $\alpha = 1$ with that evaluated at $\alpha =\hat \lambda$ and $\lambda =1$. We then observe that
$$
\textup{I}(\hat x, \mathbb{S}^3, \lambda^2 g_0) = \hat\lambda^{-k} \textup{I}(\hat x, \mathbb{S}^3,  g_0).
$$
Therefore it is a conformal scalar invariant for the surface $\hat x$ in the 3-sphere.
\endproof

\begin{appendix}

\section{The inverse of $I^{\tilde x}$ in general parametrizations}\label{app-inverse-tilde}

We consider the general parametrization 
$$
\tilde x = \alpha y_\lambda + \alpha \rho y^*_\lambda: \mathbb{R}^+\times\mathbb{R}^+\times\textup{M}^2\to \mathbb{R}^{1, 4}.
$$
Then the first fundamental form in matrix form is 
\begin{equation}\label{matrix-tilde-app}
I_{\tilde x} = \left[\begin{matrix} \begin{matrix} \ -2\rho & -\alpha \\ -\alpha & \ 0\end{matrix} & \begin{matrix} 0 & 0 \\
\alpha^2\omega^\lambda_1 & \alpha^2 \omega^\lambda_2\end{matrix}\\
\begin{matrix} \quad\quad 0 & \quad \alpha^2\omega^\lambda_1\\ \quad\quad 0 & \quad \alpha^2\omega^\lambda_2\end{matrix} & \quad 
\alpha^2 F \quad\end{matrix}\right]
 \end{equation}
 where
 \begin{equation}\label{F-matrix-app}
 \left\{
 \aligned
 F_{11} & = \frac 1m(p^2 + q^2)  +2 \rho (\omega^\lambda_1)^2\\
 F_{12} & =  F_{21} =  \frac 1m q(p+r)  + 2 \rho \omega^\lambda_1\omega^\lambda_2\\
  F_{22} & = \frac 1m(q^2 + r^2)  + 2 \rho(\omega^\lambda_2)^2 \endaligned\right. \text{ and } 
  \left\{
 \aligned
 F^*_{11} & = \frac 1m(p^2 + q^2) \\
 F^*_{12} & =  F_{21} =  \frac 1m q(p+r)  \\
  F^*_{22} & = \frac 1m(q^2 + r^2)   \endaligned\right. 
\end{equation}
 and
$$
\left[\begin{matrix} p & q \\ q & r\end{matrix}\right] = \Omega_\lambda + \rho\Omega^*_\lambda.
$$
It is easily seen that 
\begin{equation}\label{F-star-inverse}
(F^*)^{-1} = \frac m {(pr-q^2)^2}\left[\begin{matrix} r^2 + q^2 & - q(p+r)\\-q(p+r) & p^2 + q^2\end{matrix}\right]
\end{equation}
and
$$
F|_{\rho =0} = F^*|_{\rho=0} = E\left[\begin{matrix} 1 & 0 \\0 & 1\end{matrix}\right].
$$
Let 
$$
(I^{\tilde x})^{-1} = \left[\begin{matrix} a_{11} & a_{12} & a_{13} & a_{14} \\
a_{21} & a_{22} & a_{23} & a_{24} \\ a_{31} & a_{32} & a_{33} & a_{34} \\
a_{41} & a_{42} & a_{43} & a_{44} \end{matrix}\right].
$$
Therefore, for example, 
\begin{equation}\label{inverse-1-app}
 \left\{\aligned
 -2\rho a_{11} -\alpha a_{12} \quad \quad & = 1 \\
 -\alpha a_{11} + \alpha^2\omega_1a_{13} + \alpha^2\omega_2a_{14} & = 0\\
 \alpha^2\omega_1a_{12} + \alpha^2F_{11}a_{13} + \alpha^2F_{21}a_{14} & = 0 \\
 \alpha^2\omega_2a_{12} + \alpha^2F_{12}a_{13} + \alpha^2F_{22}a_{14} & = 0.
 \endaligned\right.
 \end{equation}
Subtracting the first equation multiplied by $\alpha$ from the second equation multiplied by 2 in \eqref{inverse-1-app},  we get
\begin{equation}\label{5-equ-app}
 \alpha^2 a_{12} + 2\alpha^2\rho\omega_1a_{13} + 2\alpha^2\rho\omega_2a_{14} = -\alpha
\end{equation}
And subtracting \eqref{5-equ-app} multiplied with $\omega_1$ from the third equation in \eqref{inverse-1-app} as well as subtracting \eqref{5-equ-app}
multiplied with $\omega_2$ from the fourth equation in \eqref{inverse-1-app}, we get
\begin{equation}\label{app-13-14}
\alpha^2F^*\left[\aligned a_{13}\\a_{14}\endaligned\right] = \left[\aligned \alpha\omega_1\\ \alpha\omega_2\endaligned\right]
\end{equation}
Plugging back what are $a_{13}$ and $a_{14}$ to the equation \eqref{5-equ-app} we have
\begin{equation}\label{app-12}
\left\{\aligned
a_{12}  & =  \alpha^{-1}( -1   - 2\rho [\omega_1,  \omega_2](F^*)^{-1}\left[\begin{matrix} \omega_1\\ \omega_2\end{matrix}\right])\\
a_{11}  &  = - \frac { \alpha a_{12} +1} {2\rho} = [\omega_1,  \omega_2](F^*)^{-1}\left[\begin{matrix} \omega_1\\ \omega_2\end{matrix}\right] .
\endaligned\right.
\end{equation}
Similarly one gets
\begin{equation}\label{app-2j}
\alpha^2F^* \left[\aligned a_{23}\\a_{24}\endaligned\right] = \left[\aligned-2\rho\omega_1\\ -2\rho\omega_2\endaligned\right]
\text{ and } \left\{\aligned
a_{21}  & =  \alpha^{-1}( -1   - 2\rho [\omega_1,  \omega_2](F^*)^{-1}\left[\begin{matrix} \omega_1\\ \omega_2\end{matrix}\right])\\
a_{22} & =  \frac {2\rho}{\alpha^2} ( 1 +2\rho [\omega_1,  \omega_2](F^*)^{-1}\left[\begin{matrix} \omega_1\\ \omega_2\end{matrix}\right])
\endaligned\right.
\end{equation}
\begin{equation}\label{app-3j}
\alpha^2F^* \left[\aligned a_{33}\\a_{34}\endaligned\right] = \left[\aligned 1 \\ 0 \endaligned\right] \text{ and } \left\{
\aligned
a_{31}  & =  \alpha^{-1} [\omega_1,  \omega_2](F^*)^{-1}\left[\begin{matrix} 1\\ 0\end{matrix}\right]\\
a_{32} & =  - \frac {2\rho}{\alpha^2} ( [\omega_1,  \omega_2](F^*)^{-1}\left[\begin{matrix}1 \\ 0 \end{matrix}\right]
\endaligned\right.
\end{equation}
\begin{equation}\label{app-4j}
\alpha^2F^* \left[\aligned a_{43}\\a_{44}\endaligned\right] = \left[\aligned 0\\ 1\endaligned\right]\text{ and } \left\{
\aligned
a_{41}  & =  \alpha^{-1}(F^*)^{-1}\left[\begin{matrix} 0 \\1\end{matrix}\right] \\
a_{42} & =  \frac {2\rho}{\alpha^2} ( 1 +2\rho [\omega_1,  \omega_2](F^*)^{-1}\left[\begin{matrix} \omega_1\\ \omega_2\end{matrix}\right])
\endaligned\right.
\end{equation}

\section{The geometry of the 3-sphere $\mathbb{S}^3_\lambda$ in $\mathbb{R}^{1,4}$}\label{gauss}

Let us calculate the Gauss Theorem for the 3-sphere $\mathbb{S}^3_\lambda$ in Minkowski spacetime $\mathbb{R}^{1,4}$. There is nothing new or difficult about the calculation,
but this helps to understand better about the geometry of the 3-sphere $\mathbb{S}^3_\lambda\subset\mathbb{N}^4_+\subset\mathbb{R}^{1,4}$. It is very crucial and important in our approach to use the 
fact that the induced metric on $\mathbb{S}^3_\lambda$ is exactly the conformal metric $\lambda^2g_0$. We consider the Fermi parametrization induced from 
a parametrization of the surface $\hat x: \textup{M}^2\to \mathbb{S}^3$ such that 
\begin{equation}\label{lambda-3-sphere}
y_\lambda = \lambda (\hat x(u^1, u^2, u^3))(1, \hat x(u^1, u^2, u^3)): \textup{M}^3\to \mathbb{S}^3_\lambda\subset\mathbb{N}^4_+\subset\mathbb{R}^{1,4}
\end{equation}
with 
\begin{equation}\label{extension-y}
\hat x(u^1, u^2, 0) = \hat x(u^1, u^2) \text{ and } (y_\lambda)_{u^3}|_{u^3=0} = \overset{\rightarrow}{\bf n}_\lambda. 
\end{equation}
Notice that $y_\lambda$ here is the extension of  $\hat\lambda(1, \hat x)$ before. We use the two null normal vectors $\{y_\lambda, y_\lambda^\dagger\}$ where
\begin{equation}\label{y-dagger}
<y^\dagger_\lambda, y_\lambda> -1, <y^\dagger_\lambda, (y_\lambda)_{u^1}> = <y^\dagger_\lambda, (y_\lambda)_{u^2}>  = <y^\dagger_\lambda, (y_\lambda)_{u^3}> =0.
\end{equation}
The first fundamental form is 
\begin{equation}\label{lambda-sphere-I}
I^{\mathbb{S}^3_\lambda} = \lambda^2g_0 = <dy_\lambda, dy_\lambda>.
\end{equation}
And the second fundamental form is
\begin{equation}\label{lambda-sphere-II}
II^{\mathbb{S}^3_\lambda} = - <dy_\lambda, dy^\dagger_\lambda> y^\dagger_\lambda - <dy_\lambda, dy_\lambda>y_\lambda
\end{equation}
To calculate the curvature for the metric $g_\lambda = \lambda^2g_0$ we calculate
$$
\nabla^\lambda_{\partial_{u^j}}\nabla^\lambda_{\partial_{u^i}}\partial_{u^k} - \nabla^\lambda_{\partial_{u^i}}\nabla^\lambda_{\partial_{u^j}}\partial_{u^k} 
= R^\lambda(\partial_{u^i}, \partial_{u^j})\partial_{u^k} = (R^\lambda)_{ijk}^{\quad \ l}\partial_{u^l}.
$$
First
$$
\nabla^\lambda_{\partial_{u^j}}\partial_{u^k} = (y_\lambda)_{u^ku^j} - <(y_\lambda)_{u^j}, (y^\dagger_\lambda)_{u^k}>y_\lambda - <(y_\lambda)_{u^j}, (y_\lambda)_{u^k}>y^\dagger_\lambda
$$ 
Then
$$
\aligned
\partial_{u^i}\nabla^\lambda_{\partial_{u^j}}\partial_{u^k} & = (y_\lambda)_{u^ku^ju^i} - <(y_\lambda)_{u^j}, (y^\dagger_\lambda)_{u^k}>_{u^i} y_\lambda 
- <(y_\lambda)_{u^j}, (y_\lambda)_{u^k}>_{u^i}y^\dagger_\lambda \\
& \quad -<(y_\lambda)_{u^j}, (y^\dagger_\lambda)_{u^k}>(y_\lambda)_{u^i}  - <(y_\lambda)_{u^j}, (y_\lambda)_{u^k}>(y^\dagger_\lambda)_{u^i} 
\endaligned
$$
and
$$
\aligned
\nabla^\lambda_{\partial_{u^i}}\nabla^\lambda_{\partial_{u^j}}\partial_{u^k}  & = (\partial_{u^i}\nabla^\lambda_{\partial_{u^j}}\partial_{u^k})^{T\mathbb{S}^3_\lambda}\\
& = (y_\lambda)_{u^ku^ju^i}^{T\mathbb{S}^3_\lambda}  -<(y_\lambda)_{u^j}, (y^\dagger_\lambda)_{u^k}>(y_\lambda)_{u^i}  - <(y_\lambda)_{u^j}, (y_\lambda)_{u^k}>(y^\dagger_\lambda)_{u^i} 
\endaligned
$$
Hence
$$
\aligned
(R^\lambda)_{ijk}^{\quad \ l}\partial_{u^l} & =  <(y_\lambda)_{u^j}, (y^\dagger_\lambda)_{u^k}>(y_\lambda)_{u^i}  +  <(y_\lambda)_{u^j}, (y_\lambda)_{u^k}>(y^\dagger_\lambda)_{u^i} \\
& \quad  -<(y_\lambda)_{u^i}, (y^\dagger_\lambda)_{u^k}>(y_\lambda)_{u^j}  - <(y_\lambda)_{u^i}, (y_\lambda)_{u^k}>(y^\dagger_\lambda)_{u^j} \\
\endaligned
$$
One may realize that
$$
<(y^\dagger_\lambda)_{u^i}, y^\dagger_\lambda > = 0 \text{ and } <(y^\dagger_\lambda)_{u^i}, y_\lambda>  = 0
$$
and conclude
$$
(y^\dagger_\lambda)_{u^i} = (g_\lambda)^{ml}<(y^\dagger_\lambda)_{u^i}, (y_\lambda)_{u^m}>(y_\lambda)_{u^l}.
$$
Therefore
$$
\aligned
(R^\lambda)_{ijk}^{\quad \ l} \partial_{u^l} & =  (<(y_\lambda)_{u^j}, (y^\dagger_\lambda)_{u^k}>\delta_i^{\ l}  +  (g_\lambda)_{jk}(g_\lambda)^{ml}<(y^\dagger_\lambda)_{u^i}, (y_\lambda)_{u^m}>\\
 &\quad - <(y_\lambda)_{u^i}, (y^\dagger_\lambda)_{u^k}>\delta_j^{\ l}  -  (g_\lambda)_{ik}(g_\lambda)^{ml}<(y^\dagger_\lambda)_{u^j}, (y_\lambda)_{u^m}>)
 \partial_{u^l}
 \endaligned
 $$
and
$$
\aligned
(R^\lambda)_{ijkl} = & (R^\lambda)_{ijk}^{\quad \ n}(g_\lambda)_{nl} =  <(y_\lambda)_{u^j}, (y^\dagger_\lambda)_{u^k}>(g_\lambda)_{il} 
+ <(y^\dagger_\lambda)_{u^i}, (y_\lambda)_{u^l}> (g_\lambda)_{jk}\\
 &\quad - <(y_\lambda)_{u^i}, (y^\dagger_\lambda)_{u^k}> (g_\lambda)_{jl}  -  <(y^\dagger_\lambda)_{u^j}, (y_\lambda)_{u^l}> (g_\lambda)_{ik}.
 \endaligned
 $$
 On the surface $\hat x$,  where $u^3=0$, we have
$$
 [(g_\lambda)_{ij}] = \left[\begin{matrix} E_\lambda & 0 & 0 \\0 & E_\lambda & 0\\0 & 0 & 1\end{matrix}\right].
$$
Therefore we have, for $i, j\in\{1, 2\}$,
$$
\left\{\aligned - < (y_\lambda)_{u^i}, (y_\lambda^\dagger)_{u^j}> - < (y_\lambda)_{u^3}, (y_\lambda^\dagger)_{u^3}>E_\lambda \delta_{ij} & = (R^\lambda)_{i3j3}\\
- <(y_\lambda)_{u^j}, (y^\dagger_\lambda)_{u^3}> E_\lambda\delta_{jl} + < (y_\lambda)_{u^l}, (y_\lambda^\dagger)_{u^3}> E_\lambda & = (R^\lambda)_{3jjl}\\
- < (y_\lambda)_{u^i}, (y_\lambda^\dagger)_{u^i}> E_\lambda - < (y_\lambda)_{u^j}, (y_\lambda^\dagger)_{u^j}> E_\lambda & = (R^\lambda)_{ijij} \endaligned
\right.
$$
Finally we obtain, for $i, j\in\{1, 2\}$,
\begin{equation}\label{n-dagger}
<\overset{\rightarrow}{\bf n}_\lambda, (y^\dagger_\lambda)_{u^i}> = \frac 1{E_\lambda} (R^\lambda)_{ijj3} = - (R^\lambda)_{i3},
\end{equation}
and for $i\neq j$, 
\begin{equation}\label{i-dagger-j}
\aligned
<(y_\lambda)_{u^i}, (y^\dagger_\lambda)_{u^j}>  &= - (R^\lambda)_{i3j3}\\
 <(y_\lambda)_{u^i}, (y^\dagger_\lambda)_{u^i}> & =  - (R^\lambda)_{i3i3} + \frac 12((R^\lambda)_{33} - (R^\lambda)_{1212})\\
 <(y_\lambda)_{u^3}, (y^\dagger_\lambda)_{u^3}> & = - \frac 12((R^\lambda)_{33} - (R^\lambda)_{1212})\endaligned
 \end{equation}
 
Finally, for the induced Fermi coordinate from an isothermal coordinate, we can easily see that
\begin{equation}\label{coord-covar}
\aligned
(R^\lambda)_{3i,}^{\quad i} & =  \frac 1{E_\lambda}(\sum_{i=1}^2R^\lambda)_{3i, i} \\
& =  \frac 1{E_\lambda}\sum_{i=1}^2(((R^\lambda)_{3i})_{u^i} - (R^\lambda)_{3k}(\Gamma_\lambda)^k_{ii}) \\
& =  \frac 1{E_\lambda}\sum_{i=1}^2(((R^\lambda)_{3i})_{u^i}
\endaligned
\end{equation}
Because $\sum_{i=1}^2(\Gamma_\lambda)^k_{ii} = 0$ for each $k= 1,2$, where $(\Gamma_\lambda)^k_{ij}$ is the Christofel symbols 
for the conformal metric $I_\lambda = E_\lambda |du|^2$ in the isothermal coordinates.
\end{appendix}

\end{document}